\documentclass[10pt]{article}

\usepackage{amsmath, amssymb, paralist}
\usepackage{caption,graphicx}
\usepackage[text={6.75in,9.5in},centering,letterpaper]{geometry}
\usepackage[comma,square,sort&compress,numbers]{natbib}
\usepackage{subcaption}
\usepackage{csquotes}

\setcaptionmargin{0.25in}

\newtheorem{Lemma}{Lemma}[section]
\newtheorem{Corollary}[Lemma]{Corollary}

\newtheorem{Hypothesis}{Hypothesis}
\newtheorem{Proposition}[Lemma]{Proposition}
\newtheorem{Remark}[Lemma]{Remark}
\newtheorem{Theorem}{Theorem}

\newenvironment{Proof}[1][.]%
 {\begin{trivlist}\item[]\textbf{Proof#1 }}%
 {\hspace*{\fill}$\rule{0.3\baselineskip}{0.35\baselineskip}$\end{trivlist}}

\newenvironment{Acknowledgment}%
 {\begin{trivlist}\item[]\textbf{Acknowledgments }}{\end{trivlist}}

\makeatletter\@addtoreset{equation}{section}\makeatother

\def\Re{\mathop\mathrm{Re}\nolimits}    
\def\Im{\mathop\mathrm{Im}\nolimits}    

\newcommand{\C}{\mathbb{C}}             
\newcommand{\R}{\mathbb{R}}             
\newcommand{\Z}{\mathbb{Z}}             

\newcommand{\rmO}{\mathrm{O}}           
\newcommand{\rmo}{\mathrm{o}}           
\newcommand{\rmd}{\mathrm{d}}           
\newcommand{\rme}{\mathrm{e}}           
\newcommand{\rmi}{\mathrm{i}}           
\newcommand{\per}{_\mathrm{per}}           
\newcommand{\s}{\mathrm{s}}           
\newcommand{\un}{\mathrm{u}}           
\newcommand{\opno}{_{L(X^0,X^0\per)}}   
\newcommand{\xper}{_{X^0\per}}	

\begin{document}

\title{Behavior of Spiral Wave Spectra with a Rank-Deficient Diffusion Matrix}

\author{%
Stephanie Dodson\\
Department of Mathematics\\
University of California, Davis\\
Davis, CA~95616, USA
\and
Bj\"orn Sandstede\\
Division of Applied Mathematics\\
Brown University\\
Providence, RI~02912, USA
}

\date{\today}
\maketitle

\begin{abstract}
Spiral waves emerge in numerous pattern forming systems and are commonly modeled with reaction-diffusion systems. Some systems used to model biological processes, such as ion-channel models, fall under the reaction-diffusion category and often have one or more non-diffusing species which results in a rank-deficient diffusion matrix. Previous theoretical research focused on spiral spectra for strictly positive diffusion matrices. In this paper, we use a general two-variable reaction-diffusion system to compare the essential and absolute spectra of spiral waves for strictly positive and rank-deficient diffusion matrices. We show that the essential spectrum is not continuous in the limit of vanishing diffusion in one component. Moreover, we predict locations for the absolute spectrum in the case of a non-diffusing slow variable. Predictions are confirmed numerically for the Barkley and Karma models. 
\end{abstract}

\section{Introduction}
 
Spiral waves are frequently observed in nature, including chemical oscillations in the Belousouv-Zhabotinsky reaction \cite{Zaikin:1970te,Winfree:1972ud}, cell signaling patterns in slime molds \cite{Newell:1982jr}, and in electrical activity in cardiac dynamics \cite{Winfree:1994cz,Rosenbaum:1994bt}. Stable spiral waves are observed in these systems, but bifurcations to complex and unstable pattens are common and these bifurcations can have profound results. For example, rotating spiral waves in cardiac electrical activity have been linked to dangerous tachycardiac rhythms and the transition to break up can lead to life-threatening fibrillation and sudden cardiac death \cite{Rosenbaum:1994bt,Pastore:1999dq}. Therefore, understanding the stability and bifurcations of spiral waves poses interesting mathematical questions and is important in applications.

Reaction-diffusion systems are canonical pattern forming systems describing biological and physical processes and take on the form
\begin{align}\label{eqn:rxn_diff_eqn}
U_t = D \Delta U + F(U), \qquad
U = U(y,t) \in \mathbb{R}^n, \ D \in \mathbb{R}^{n \times n}, \ y \in \mathbb{R}^2,
\end{align}
where $\Delta$ is the Laplacian and the smooth, typically nonlinear $F(U)$ defines kinetic reaction terms. The $n$ species of $U = [u_1, \dots, u_n]^T$ diffuse with diffusion rates given by the entries $D_{jj}$ of the positive diagonal matrix $D$ for $j=1,\ldots,n$. Planar spiral waves have a regular shape and rotate with angular frequency $\omega_0$. Thus, these pattens are stationary solutions of (\ref{eqn:rxn_diff_eqn}) in a polar coordinate rotating frame, and stability can be investigated by evaluating the spectrum of the reaction-diffusion operator linearized about the spiral wave solution. 

The propagation of electrical potentials in excitable media, such as neurons and cardiac tissue, can be described by biophysically detailed ion-channel models given by the system
\[
V_t = \Delta V + f(V,n), \qquad
n_t = g(V,n). 
\]
Here, $V = V(y,t) \in \mathbb{R}$ corresponds to the electrical potential and $n = \left(n_1(y,t),\dots, n_M(y,t) \right)^T\in~\mathbb{R}^M$ are $M$ dynamic gating variables explaining the opening and closing of ion channels that facilitate the voltage propagation. Ion channel models can still be written in the general reaction-diffusion framework, with the popular Hodgkin--Huxley \cite{Hodgkin:1945ef}, Noble \cite{Noble:1962tk}, and Beeler--Reuter \cite{Beeler:1977ct} models falling into this category.

\begin{figure}
\centering
 \includegraphics[width=0.75\linewidth]{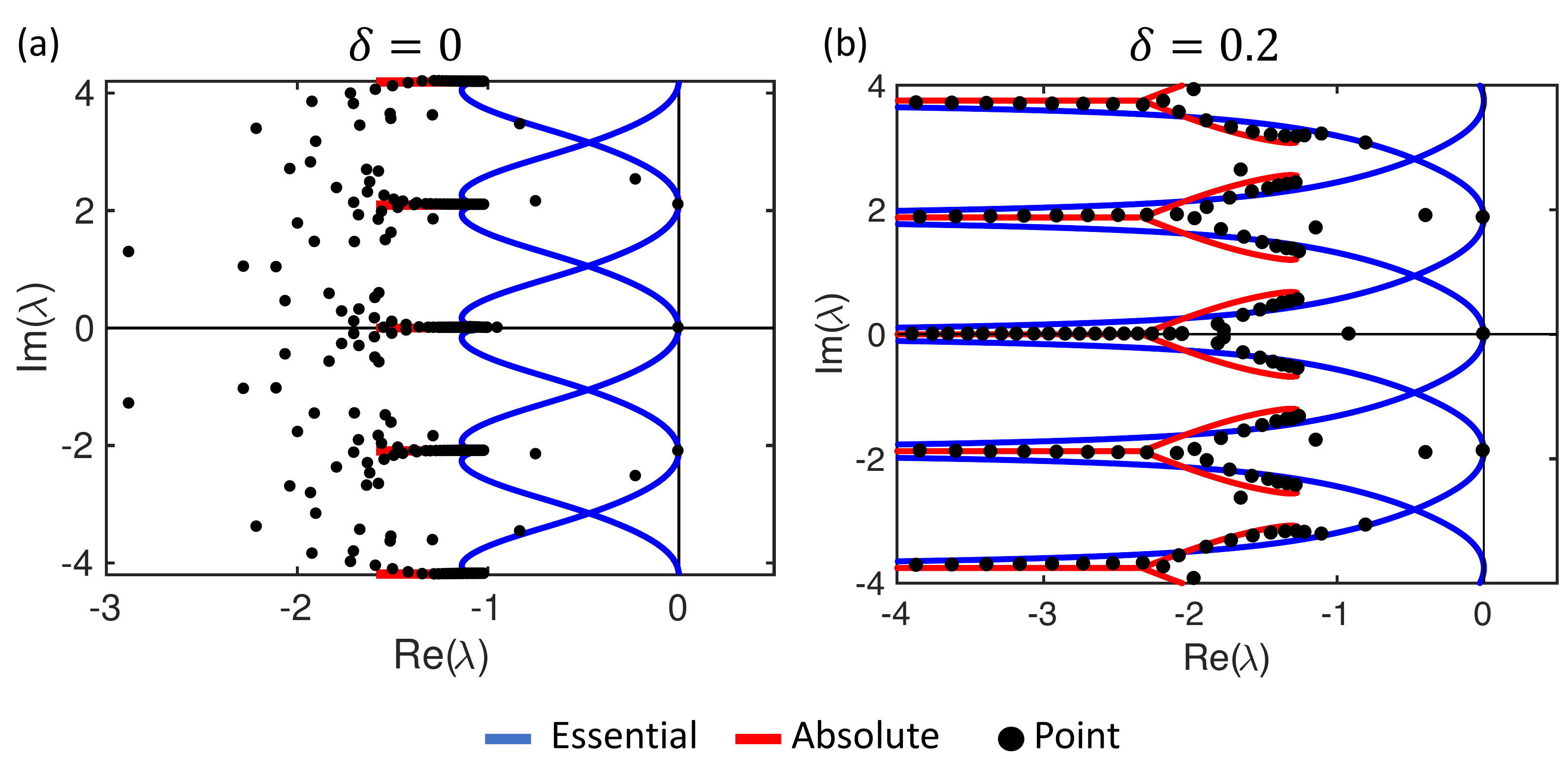}
\footnotesize \caption{Barkley Model: Differences in the essential, absolute, and point spectra for $\delta>0$ and $\delta= 0$.} \label{fig:fig1}
\end{figure}

In this paper, we focus on two-component reaction-diffusion models. In ion channel models, the gating variable does not diffuse and an appropriate ion channel models is therefore of the form
\[
u_t = \Delta u + f(u,v), \qquad
v_t = g(u,v),
\]
where we now use $(u,v)$ instead of $(V,n)$. Often, small unphysical diffusion is added to these components or included in qualitative models. Throughout, we will use $\delta$ to correspond to the small diffusion coefficient of interest, leading to the adjusted model
\[
u_t = \Delta u + f(u,v), \qquad
v_t = \delta\Delta v + g(u,v)
\]
where $0<\delta\ll1$. We would expect that the spectra should change smoothly as the diffusion coefficient $\delta \rightarrow 0$, and this is indeed true for the continuous spectra of one-dimensional periodic wave trains \cite{Rademacher:2004ue}. Yet, we observe that spiral spectra computed for $\delta>0$ do not converge as $\delta \searrow 0$ to the spectrum computed for $\delta=0$. Figure~\ref{fig:fig1} shows the differences in the spectra for the Barkley model, a two-variable reaction-diffusion system, with $\delta = 0.2$ and $\delta = 0$. For $\delta>0$, the spectral curves are unbounded, whereas the curves remain bounded for $\delta = 0$. Furthermore at the finite limit points, adjacent curves meet and form the cusps seen in Figure~\ref{fig:fig1}.  Additional changes are seen in the absolute and point spectrum, with the absolute spectrum collapsing to short branches that align with the essential spectrum cusp points. The effect of a rank-deficient diffusion matrix on the structure of the spectra of planar spiral waves has previously not been analyzed.

In this paper, we focus on changes in the continuous spectrum for spiral waves in the $\delta=0$ limit and investigate the mechanisms responsible for the discontinuities of essential spectrum observed in the limit of a rank-deficient diffusion matrix. Under mild conditions on the linearization of the reaction terms, Theorems~\ref{t:1}-\ref{t:2} provide expressions for the spectral curves near the limit points. Our findings indicate the discontinuity occuring in the $\delta = 0$ limit is dictated by the non-diffusing species. Furthermore, we use the result of Theorem~\ref{t:1} to predict locations for the rank-deficient absolute spectrum in Theorem~\ref{thm:abs_spec}.

We proceed as follows:  Mathematical preliminaries of spiral waves and the relevant spectral properties are described in Section~\ref{sec:math}. The main results for the continuous spectra of planar spiral waves are derived in  Sections~\ref{s:3}-\ref{s:4}.  Implications for the locations of absolute spectrum are discussed in Section~\ref{sec:abs_spec}. Finally, the results are applied to the Barkley and Karma models in Section \ref{sec:barkley_karma}. 

\section{Mathematical preliminaries} \label{sec:math}

\subsection{Planar spiral waves and asymptotic wave trains}
Planar spiral waves are solutions to the reaction-diffusion equation~(\ref{eqn:rxn_diff_eqn}) which rigidly rotate in time with a constant angular temporal frequency $\omega_0$. Spiral waves are therefore stationary solutions $U_*(r,\psi)$ in the co-rotating frame under polar coordinates $(r,\psi)$
\begin{align*}
U_t = D \Delta_{r,\psi} U + \omega_0 U_{\psi} + F(U), \qquad U \in \mathbb{R}^n
\end{align*}
where $\Delta_{r,\psi}$ denotes the Laplacian expressed in polar coordinates. In the limit $r \rightarrow \infty$, the spiral solutions $U_*(r,\psi)$ limit to $2\pi/\kappa$-periodic functions $U_{\infty}(r + \psi/\kappa) = U_{\infty}(x)$ referred to as asymptotic wave trains, which are stationary solutions to 
\begin{align*}
U_t = D U_{xx} + \omega U_{x} + F(U)
\end{align*}
where $\omega = \omega_0/\kappa$. We focus on planar spiral waves that can be viewed as a source which emits wave trains with positive group velocity. The spatial wave number $\kappa$ is selected by the spiral and the nonlinear dispersion relation $\omega = \omega_*(\kappa)$ of the wave train connects $\omega$ and $\kappa$. Our main assumption is the existence of a spiral-wave solution that depends smoothly on the diffusion coefficient $\delta$ for $0\leq\delta\ll1$.
\begin{Hypothesis}\label{hyp:1}
The spiral-wave solution $U_*$ and the asymptotic wave train $U_{\infty}$ both depend smoothly on the diffusion coefficient $\delta$ for $0\leq\delta\ll1$.
\end{Hypothesis}
While we do not have a proof that Hypothesis~\ref{hyp:1} holds, this hypothesis is well supported by numerical evidence.

\subsection{Essential spectrum of wave trains}
We are interested in the continuous spectrum of the reaction-diffusion operator linearized about a spiral wave pattern. The reader is referred to \cite{Fiedler:2000,Sandstede:2000ug,Sandstede:2002ht,Kapitula:2013} for a detailed study of spectrum of operators in nonlinear waves. 

We start with the wave train spectrum as it is directly linked to the continuous spectrum of the spiral wave. The spectrum of the linearization
\[
\mathcal{L}_{\infty}V =  D V_{xx} + \omega V_{x} + F_U(U_{\infty})V = \lambda_{\infty} V
\]
of the asymptotic wave trains $U_{\infty}(x)$ on the space $L^2(\mathbb{R},\mathbb{R}^n)$ is given by the set of $\lambda_{\infty} \in \mathbb{C}$ for which the eigenvalue problem $\mathcal{L}_{\infty}V = \lambda_{\infty} V$ has a nontrivial solution of the form $V(x) = \rme^{\nu x} \bar{V}(x)$ where $\bar{V}(x+ 2\pi/\kappa) = \bar{V}(x)$ and $\nu = \rmi \gamma$ for $\gamma \in \mathbb{R}$. Wave train eigenvalues are therefore determined by nontrivial $2\pi/\kappa$-periodic solutions $\bar{V}(x)$ of
\begin{align} \label{eqn:linear_disp_wt}
\mathcal{L}_{\infty}(\lambda,\nu) \bar{V} = D \left( \partial_{x} + \nu\right)^2\bar{V} + \omega \left(  \partial_{x} + \nu \right) \bar{V} + F_U(U_{\infty}) \bar{V} - \lambda_{\infty} \bar{V} = 0
\end{align}
and come in curves $\lambda_{\infty} = \lambda_{\infty}(\rmi \gamma)$ parameterized by the Floquet exponent $\nu = \rmi \gamma$. From translational symmetry, one spectral curve includes $\lambda_{\infty}(0) = 0$ for $\nu = 0$ with eigenfunction $V(x) = U'_{\infty}(x)$. The orientation of each curve is defined as the direction of increasing $\gamma$. 

\subsection{Essential spectrum of planar spiral waves}
The spectrum of a planar spiral wave is given by considering the spectrum of the operator
\begin{align} \label{eqn:spiral_eval}
\mathcal{L}_* V = D \Delta_{r,\psi} V + \omega_0 V_{\psi} + F_U(U_*)V
\end{align}
on $L^2(\mathbb{R},\mathbb{R}^n)$ which contains isolated eigenvalues in the point spectrum and continuous curves of essential spectrum. The essential spectrum of the spiral wave is determined by the far-field asymptotic dynamics, and it can be shown that the spiral wave continuous spectral curves $\lambda_*(\nu)$ are related to those of the wave trains via 
\begin{align} \label{eqn:wt_spiral_relationship}
\lambda_*(\nu) = \lambda_{\infty}(\nu) - \frac{\omega_0}{\kappa} \nu + \rmi \omega_0 \ell, \ \ \ell \in \mathbb{Z}.
\end{align}
Therefore, the linear dispersion relation for the spiral wave is determined by the solvability of
\begin{align} \label{eqn:sp_disp_rel}
\mathcal{L}_{*}(\lambda,\nu) \bar{V} = D \left( \partial_x + \nu\right)^2 \bar{V}+ \omega  \bar{V}_{x} + F_U(U_{\infty}) \bar{V} - \lambda_* \bar{V} = 0.
\end{align}
Since $\nu \in \rmi \mathbb{R}$, the real parts wave train and spiral wave continuous spectra coincide. The $\rmi \omega_0 \ell$ term provides additional vertically periodic branches due to the rotational symmetry of the spiral wave. 

\subsection{Computation of essential spectra}
In the remainder of the paper, we focus on spiral waves described by a two-component reaction-diffusion system $U = (u,v)^T$ set in a rotating frame with polar coordinates $(r,\psi) = (r,\phi - \omega t)$, 
\begin{align} \label{eqn:2var_rxn_diff}
u_t = & \ \Delta_{r,\psi} u + \omega u_{\psi} + f(u,v)\\
v_t = & \ \delta \Delta_{r,\psi} v + \omega v_{\psi} + g(u,v) \nonumber
\end{align}
where $0 \leq \delta \ll 1$. In this framework, $u$ and $v$ are typically referred to as the fast and slow diffusing species, respectively. Many of the reduced systems, including the Barkley, Karma, and Morris-Lecar models, have two components.

In these systems, the equation for the essential spectrum of the planar spiral wave is then given by
\begin{align} \label{eqn:2comp_rxn_diff_ess_spec}
\lambda_* u =  &\left(  \partial_{x} + \rmi\gamma \right)^2 u + \omega u_{x} + f_u\left(U_{\infty}(x) \right)u + f_v\left(U_{\infty}(x) \right)v\\
\lambda_* v = \delta &\left(  \partial_{x} + \rmi\gamma \right)^2 v + \omega v_{x} + g_u\left(U_{\infty}(x) \right)u + g_v\left(U_{\infty}(x) \right)v, \nonumber
\end{align} 
where the imaginary part $\gamma$ of the spatial Floquet exponent parameterizes the essential spectrum curves, and the nonlinear terms are linearized about the wave train $U_{\infty}(x)$. In Theorems~\ref{t:1}-\ref{t:2}, we assume that the derivative $g_v\left(U_{\infty}(x) \right) = \bar{g}$ is a constant. To simplify notation, we remove the explicit $U_{\infty}$ dependence in the linearizations and instead write $f_{u,v}\left(U_{\infty}(x) \right) = f_{1,2}(x)$ and $g_{u}\left(U_{\infty}(x) \right) = g_{1}(x)$. The observed changes in the essential spectrum are for large values of $\gamma$ so that $|\gamma| \gg 1$. Sending $\gamma \rightarrow \pm \infty$ corresponds to traversing the curve in opposite directions. To focus on the limiting behavior (and the area near the cusp point for $\delta = 0$),we  divide the system (\ref{eqn:2comp_rxn_diff_ess_spec}) by $\gamma^2$, define $\alpha = 1/\gamma$, and group terms by orders of $\alpha$ to arrive at
\begin{align} \label{eqn:dispersion_relation_alpha_u}
-u + \alpha \left( 2\rmi u_{x} \right) + \alpha^2 \left(  u_{xx} + \omega u_{x} + f_1(x)u  + f_2(x) v- \lambda u \right) = 0\\
-\delta v + \alpha \left( 2 \rmi \delta  v_{x} \right) + \alpha^2 \left( \delta  v_{xx} + \omega v_{x} + g_1(x)u  + \bar{g} v - \lambda v \right) = 0  \label{eqn:dispersion_relation_alpha_v}
\end{align}
where $(u,v)(x + 2\pi/\kappa) = (u,v)(x)$. The theorems and lemmas to follow analyze the continuous spectrum assuming the conditions in Hypothesis~\ref{hyp:1} and using the formulation in equations (\ref{eqn:dispersion_relation_alpha_u})-(\ref{eqn:dispersion_relation_alpha_v}). We consider values of $|\alpha|\ll1$, where $\alpha>0$ ($\alpha <0$) corresponds to the positively (negatively) oriented tails of essential spectrum curves.

\section{Statement of main results: essential spectrum} \label{s:3}

Throughout, we will make use of \citep[Corollary~3.4]{Sandstede:2000ut}, which states that the spectrum of a planar spiral wave is invariant under the shift $\lambda\mapsto\lambda+\rmi\omega_0 n$ for integers $n\in\Z$, and therefore restrict our analysis to rectangles of the form
\begin{equation}\label{d:1}
\Lambda_R := \left\{\lambda\in\C:\; |\Im\lambda|\leq\frac12,\, |\Re\lambda|<R\right\}
\end{equation}
where $R$ is positive or $R=\infty$. 

Our first result theorem focuses on the case of vanishing diffusion in the $v$-component. For each $|\alpha|\ll1$, we seek to find all $\lambda\in\Lambda_\infty$ for which the eigenvalue problem
\begin{align}
\lambda u & = u_{xx} + \left(\frac{2\rmi}{\alpha}+\omega\right) u_x - \frac{1}{\alpha^2} u + f_1(x) u + f_2(x) v
\label{evp:1} \\ \label{evp:2}
\lambda v & = \omega v_x + \bar{g} v + g_1(x) u
\end{align}
has a nontrivial $2\pi/\kappa$-periodic complex-valued solution $(u,v)(x)$.

\begin{Theorem}\label{t:1}
Fix $R>0$, $\omega>0$, and $\bar{g}\in\R$, and assume that $f_1,f_2,g_1$ are given $2\pi/\kappa$-periodic functions of class $C^4$, then there is a constant $\alpha_0>0$ so that the following is true. There is a unique function $\lambda_0^*:[-\alpha_0,\alpha_0]\to\Lambda_R$ so that (\ref{evp:1})-(\ref{evp:2}) has a nontrivial $2\pi/\kappa$-periodic solution for $\lambda\in\Lambda_R$ and $|\alpha|\leq\alpha_0$ if and only if $\lambda=\lambda_0^*(\alpha)$. Furthermore, the function $\lambda_0^*$ is continuous, and there are constants $\lambda_{2,3}\in\R$ so that
\begin{align*}
\lambda_0^*(\alpha) = \bar{g} + \lambda_2\alpha^2 + 2\rmi\lambda_3\alpha^3 + \rmo(\alpha^3).
\end{align*}
which shows that the trace of this function in the complex plane is a cusp emerging from $\lambda_0^*(0)=\bar{g}$.
\end{Theorem}

We will provide more details of the expansions of $\lambda_0^*(\alpha)$ and the associated eigenfunctions of (\ref{evp:1})-(\ref{evp:2}) in Proposition~\ref{p:2} below. We collect two extensions of our results in the following remark.

\begin{Remark} \label{r:1}
\begin{compactitem}
\item The spectrum is invariant under the shifts $\lambda \mapsto \lambda + \rmi \omega_0 n,\  n \in \mathbb{Z}$ \cite{Sandstede:2000ut}. Thus, more generally, the expansion of $\lambda_0^*(\alpha)$ is
\begin{align*} \lambda_0^*(\alpha) = \bar{g} + \rmi \omega_0 n + \lambda_2 \alpha^2 + \rmi \lambda_3 \alpha^3 + \rmo(\alpha^4), \qquad n \in \mathbb{Z}.\end{align*}
In particular, the spectrum of a planar spiral wave with vanishing diffusion features infinitely many cusps that emerge from each $\lambda_0 = \bar{g} + \rmi \omega_0 n$ for $n \in \mathbb{Z}$. The cusp structure at each $\lambda_0$ is formed from the positive and negative tails of two curves, corresponding to $\alpha >0$ and $\alpha <0$, meeting at $\lambda_0$.
\item Theorem~\ref{t:1} also holds if we replace $\bar{g}$ by a $2\pi/\kappa$-periodic function $g_2(x)$. The expressions for the coefficients in the expansion of $\lambda_0^*(\alpha)$ remain valid provided we use the definition $\bar{g}:=\frac{\kappa}{2\pi}\int_0^{2\pi/\kappa}g_2(y)\,\rmd y$.
\end{compactitem}
\end{Remark}

Our second result focuses on the case where the equation for $v$ includes diffusion with a small diffusion constant $0<\delta\ll1$. For each fixed $0<\delta\ll1$, we will identify a range of values of $\alpha$ with $|\alpha|\ll1$ for which the eigenvalue problem
\begin{align}
\lambda u & = u_{xx} + \left(\frac{2\rmi}{\alpha}+\omega\right) u_x - \frac{1}{\alpha^2} u + f_1(x) u + f_2(x) v
\label{evp:3} \\ \label{evp:4}
\lambda v & = \delta v_{xx} + \left(\frac{2\rmi\delta}{\alpha}+\omega\right) v_x - \frac{\delta}{\alpha^2} v + \bar{g} v + g_1(x) u
\end{align}
has a nontrivial $2\pi/\kappa$-periodic complex-valued solution $(u,v)(x)$.

\begin{Theorem}\label{t:2}
Fix $\omega>0$ and $\bar{g}\in\R$, and assume that $f_1,f_2,g_1$ are given $2\pi/\kappa$-periodic functions of class $C^4$, then there are constants $\delta_0,s_0>0$ and functions $\alpha^*,\lambda^*:[-s_0,s_0]\times(0,\delta_0]\to\C$ with
\begin{align} 
\alpha^*(s,\delta) & = \frac12 \left( s+ \sqrt{s^2+4\sqrt{\delta}} \right) \\
\lambda^*(s,\delta) & = \left\{ \begin{array}{lcl} 
\lambda_0^*(|s|) + \rmO(\sqrt{\delta}) & \quad & s\geq 0 \\
\bar{g} - \frac{s^2}{\omega} (1+\rmO(\delta^{\frac14})) + \rmO(\sqrt{\delta}) & \quad & s\leq 0
\end{array} \right.
\end{align}
so that (\ref{evp:3})-(\ref{evp:4}) has a nontrivial $2\pi/\kappa$-periodic solution when $(\alpha,\lambda)=(\alpha^*,\lambda^*)(s,\delta)$.
\end{Theorem}

\begin{figure}
\centering\includegraphics{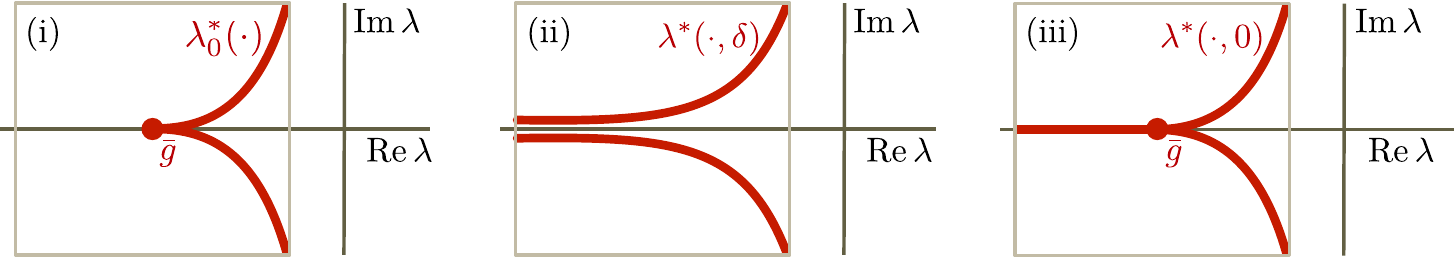}
\caption{The three panels illustrate the behavior of the dispersion curves, and their complex conjugates, of (\ref{evp:1})-(\ref{evp:2}) for $\delta=0$ and (\ref{evp:3})-(\ref{evp:4}) for fixed $0<\delta\ll1$ in a small fixed neighborhood of $\lambda=\bar{g}$. Panel~(i) shows the trace of $\lambda_0^*(\cdot)$ given in Theorem~\ref{t:1} for (\ref{evp:1})-(\ref{evp:2}), while panel~(ii) illustrates the trace of $\lambda^*(\cdot,\delta)$ established in Theorem~\ref{t:2} for (\ref{evp:3})-(\ref{evp:4}) with $\delta>0$ fixed. Panel~(iii) shows that the traces of $\lambda_0^*(\cdot)$ and $\lambda^*(\cdot,0)$ differ.}
\label{f:1}
\end{figure}

Theorem~\ref{t:2} implies that the $\delta=0$ limit of the trace of the dispersion curve $\lambda^*(\cdot,\delta)$ is given by the trace of the dispersion curve $\lambda_0^*(\cdot)$ for $\delta=0$ plus the interval $[\lambda_0^*(0)-s_0^2/\omega,\lambda_0^*(0)]$. In particular, the essential spectrum at $\delta=0$ differs from  the $\delta=0$ limit of the essential spectra for $\delta>0$ (Figure~\ref{f:1}).

\section{Proofs of main results: essential spectrum} \label{s:4}

To prove Theorems~\ref{t:1}-\ref{t:2}, we systematically construct the eigenfunction solutions to the appropriate system. Lemmas~\ref{l:1}-\ref{l:3} and Proposition~\ref{p:1} assume that $v$ is a given smooth function and establish a $2\pi/\kappa$-periodic solution $u = \mathcal{D}(\lambda,\alpha)v$ for (\ref{evp:1}). These results for the $u$-equation are independent of $\delta$, and the cases $\delta = 0$ and $\delta >0$ are then handled separately in the construction of the $v$-eigenfunction. Lemma~\ref{l:4} and Proposition~\ref{p:2} use and refine the form of $u$ to find a $2\pi/\kappa$-periodic solution for (\ref{evp:2}) when $\delta = 0$, and Theorem~\ref{t:2} follows a similar outline for the case of $\delta >0$. Throughout the proofs, we will ignore that the coefficient functions depend on $\delta$: we assumed in Hypothesis~\ref{hyp:1} that these functions depend smoothly on $\delta$ for $0\leq\delta\ll1$ and including this dependence does not change our arguments or expansions since the $\delta$-dependence shows up at higher order.

We note that (\ref{evp:3})-(\ref{evp:4}) is a singular perturbation of (\ref{evp:1})-(\ref{evp:2}) since $0\leq\delta\ll1$. However, due to the presence of the additional small parameter $0\leq\alpha\ll1$, (\ref{evp:3})-(\ref{evp:4}) cannot be rescaled so that there is a well defined limit as $(\alpha,\delta)$ approaches zero. In particular, methods such as geometric singular perturbation theory are not readily applicable.

Throughout, we will fix a constant $T>0$ and denote the stable and unstable spectral projections of a hyperbolic matrix $A$ by $P^\mathrm{s}$ and $P^\mathrm{u}$, respectively. We say a function is continuously differentiable on a closed interval if it is differentiable on the open interval, the one-sided derivatives at the boundary points exist, and the resulting derivative is continuous on the closed interval. We denote the Banach space of linear bounded operators between two given Banach spaces $X$ and $Y$ by $L(X,Y)$.

\subsection{Periodic solutions of linear systems of differential equations}

For each natural number $k\geq0$, we define the spaces
\begin{align*}
Y^k\per := C^k(\R/T\Z,\C^n), \qquad
Y^0 := C^0([0,T],\C^n), \qquad
Y^{k+1} := C^k(\R/T\Z,\C^n) \cap C^{k+1}([0,T],\C^n).
\end{align*}
Consider the inhomogeneous linear differential equation
\begin{equation}\label{e:1}
U_x = AU + G(x), \qquad U\in\C^n, \qquad 0\leq x\leq T.
\end{equation}
The following lemma establishes the existence of $T$-periodic solutions of (\ref{e:1}).

\begin{Lemma}\label{l:1}
Fix $T>0$ and assume that the matrix $A\in\C^{n\times n}$ is hyperbolic, then the following is true.
\begin{compactitem}
\item For each $G\in Y^0$ equation (\ref{e:1}) has a unique $T$-periodic solution $U=\mathcal{S}G\in Y^1$, and the operator $\mathcal{S}\in L(Y^0,Y^1)$ is given by
\begin{align*}
[\mathcal{S}G](x) =&\  \rme^{Ax}\left(1 - \rme^{AT}\right)^{-1} \int_0^T \rme^{A(T-s)} P^\s G(s) \,\rmd s\ + \int_0^x \rme^{A(x-s)} P^\s G(s) \,\rmd s\ \\
&+ \rme^{A(x-T)} \left(1 - \rme^{-AT}\right)^{-1} \int_T^0 \rme^{-As} P^\un G(s) \,\rmd s\ + \int_T^x \rme^{A(x-s)} P^\un G(s) \,\rmd s\ ,
\qquad 0\leq x\leq T.
\end{align*}
\item For each $k\geq0$, we have $\mathcal{S}\in L(Y^k\per, Y^{k+1}\per)$.
\item For each $k\geq1$, we can write $\mathcal{S}$ as
\begin{align}\label{d:s}
\mathcal{S}G = A^{-k} \mathcal{S} \frac{\rmd^k G}{\rmd x^k} - \sum_{j=1}^k A^{-j} \frac{\rmd^{j-1} G}{\rmd x^{j-1}}, \qquad
G\in Y^k\per.
\end{align}
\end{compactitem}
\end{Lemma}

\begin{Proof}
Let $E^\mathrm{s,u}:=\mathcal{R} \left(P^\mathrm{s,u}\right)$ denote the generalized stable and unstable eigenspaces of the matrix $A$. Using the variation-of-constants formula, the general solution of (\ref{e:1}) is given by
\begin{align*}
U(x) = \rme^{Ax} a^\s+ \rme^{A(x-T)} a^\un + \int_0^x \rme^{A(x-s)} P^\s G(s) \,\rmd s\ + \int_T^x \rme^{A(x-s)} P^\un G(s) \,\rmd s, \qquad 0 \leq x \leq T
\end{align*}
with $a^\mathrm{s,u}\in E^\mathrm{s,u}$. The coefficients $a^\mathrm{s,u}$ are uniquely determined upon enforcing the periodic boundary condition $U(0)=U(T)$, which yields the unique $T$-periodic solution $U(x)=[\mathcal{S}G](x)$ given by
\begin{align} \label{e:sg}
[\mathcal{S}G](x) := &\  \rme^{Ax}\left(1 - \rme^{AT}\right)^{-1} \int_0^T \rme^{A(T-s)} P^\s G(s) \,\rmd s\ + \int_0^x \rme^{A(x-s)} P^\s G(s) \,\rmd s\ \\ \nonumber
&+ \rme^{A(x-T)} \left(1 - \rme^{-AT}\right)^{-1} \int_T^0 \rme^{-As} P^\un G(s) \,\rmd s\ + \int_T^x \rme^{A(x-s)} P^\un G(s) \,\rmd s
\end{align}
for $G\in Y^0$.

Next, it follows directly from (\ref{e:1}) and periodicity of $\mathcal{S}G$ that the derivative of $\mathcal{S}G$ is periodic. In particular, $G\in Y^k\per$ implies $\mathcal{S}G\in Y^{k+1}\per$.

Finally, when $G\in Y^k\per$ for some $k\geq1$, we apply integration-by-parts $k$-times in each of the integral terms that appear in the expression for the solution operator $\mathcal{S}$. Upon each integration by parts, most of the boundary terms cancel due to the $T$-periodicity of $G$, and the remaining terms combine to provide the identity (\ref{d:s}). We omit the details of this calculation.
\end{Proof} 

\subsection{Periodic solutions of linear second-order equations}

Next, we apply these results to second-order equations of the form
\begin{equation}\label{e:2}
u_{xx} = a_{21} u + a_{22} u_x + g(x), \qquad u\in\C, \qquad 0\leq x\leq T.
\end{equation}
We write $P_1$ for the projection of $\C^2$ onto the first component and define the vector $e_2:=\begin{pmatrix} 0\\1 \end{pmatrix}\in\C^2$. Finally, we modify our function spaces and define 
\begin{align*}
X^k\per := C^k(\R/T\Z,\C), \qquad
X^0 := C^0([0,T],\C), \qquad
X^{k+1} := C^k(\R/T\Z,\C) \cap C^{k+1}([0,T],\C)
\end{align*}
for any natural number $k\geq0$.

\begin{Lemma}\label{l:2}
Fix $T>0$ and assume that $a_{21},a_{22}$ are such that the matrix $A:=\begin{pmatrix} 0 & 1 \\ a_{21} & a_{22} \end{pmatrix}\in\C^{2\times2}$ has eigenvalues $\nu^\s,\nu^\un$ with $\Re\nu^\s<0<\Re\nu^\un$, then the following is true.
\begin{compactitem}
\item For each $g\in X^0$ equation (\ref{e:2}) has a unique $T$-periodic solution $u=\mathcal{T}g\in X^2$ given by
\begin{align}
[\mathcal{T}g](x) = & \ [P_1 \mathcal{S} e_2g](x) \nonumber \\ \label{e:t}
= & \  P_1 \left[ \rme^{\nu^\s x} \left(1 - \rme^{\nu^\s T}\right)^{-1} \int_0^T \rme^{\nu^\s(T-s)} P^\s e_2 g(s) \,\rmd s\   + \int_0^x \rme^{\nu^\s (x-s)} P^\s e_2 g(s) \,\rmd s\ \right. \\ \nonumber
& \left. + \rme^{\nu^\un (x-T)} \left(1-\rme^{-\nu^\un T}\right)^{-1} \int_T^0 \rme^{-\nu^\un s} P^\un e_2 g(s) \,\rmd s\ + \int_T^x \rme^{\nu^\un (x-s)}P^\un e_2g(s) \,\rmd s\ \right], \qquad 0\leq x\leq T,
\end{align}
and we have
\begin{align}\label{e:test}
\|\mathcal{T}\|\opno \leq \frac{2}{1-\rme^{-\min(|\Re\nu^\mathrm{s,u}|)T}} \left( \frac{|P_1 P^\s e_2|}{|\Re\nu^\s|} + \frac{|P_1 P^\un e_2|}{|\Re\nu^\un|} \right).
\end{align}
\item For each $k\geq1$, we have $\mathcal{T}\in L(X^k\per,X^{k+2}\per)$ with
\begin{align}\label{d:t}
\mathcal{T}g = P_1 A^{-k} \mathcal{S} e_2 \frac{\rmd^k g}{\rmd x^k} - \sum_{j=1}^k P_1 A^{-j} e_2 \frac{\rmd^{j-1} g}{\rmd x^{j-1}}, \quad g\in X^k\per,
\end{align}
and for $g\in X^k\per$
\begin{align}\label{d:test}
\left|P_1 A^{-k} \mathcal{S} e_2 \frac{\rmd^k g}{\rmd x^k}\right|_{X^0_\mathrm{per}} \leq \frac{2}{1-\rme^{-\min(|\Re\nu^\mathrm{s,u}|)T}} \left( \frac{|P_1 P^\s e_2|}{|\Re\nu^\s| |\nu^\s|^k} + \frac{|P_1 P^\un e_2|}{|\Re\nu^\un| |\nu^\un|^k} \right) \left| \frac{\rmd^k g}{\rmd x^k}\right|\xper.
\end{align}
\end{compactitem}
\end{Lemma}

\begin{Proof}
The second-order equation (\ref{e:2}) can be rewritten as the system
\begin{align}\label{e:2s}
U_x = \begin{pmatrix} 0 & 1 \\ a_{21} & a_{22} \end{pmatrix} U + \begin{pmatrix} 0 \\ g(x) \end{pmatrix} := AU + e_2 g(x), \qquad U = \begin{pmatrix} u \\ u_x \end{pmatrix}.
\end{align}
Since $A$ is hyperbolic by assumption, Lemma~\ref{l:1} shows that (\ref{e:2s}) has a unique $T$-periodic solution in $Y^1$ for each $g\in X^0$ and that this solution is given by $U=\mathcal{S}e_2 g$. Hence, $u=P_1\mathcal{S}e_2 g=:\mathcal{T}g$ is the unique $T$-periodic solution of (\ref{e:2}) in $X^2$. Substituting $\mathcal{S}$ from Lemma~\ref{l:1} and using that $\rme^{Ax}P^\s=\rme^{\nu^\s x}P^\s$ and $\rme^{Ax}P^\un=\rme^{\nu^\un x}P^\un$ establishes (\ref{e:t}) and shows that $\mathcal{T}\in L(X^0,X^2\per)$.

To prove (\ref{e:test}), we note that
\begin{align*}
\|\rme^{Ax}|_{E^\mathrm{s}}\|\leq\rme^{\Re\nu^\mathrm{s} x} \text{ for } x \geq 0, \qquad
\|\rme^{Ax}|_{E^\mathrm{u}}\| \leq \rme^{\Re\nu^\mathrm{u} x} \text{ for } x \leq 0
\end{align*}
Hence, we obtain from (\ref{e:t}) that
\begin{align*}
| \mathcal{T} g|_{X^0} \leq &\ \max_{0 \leq x \leq T}
\left[ \left(\frac{\rme^{\Re \nu^\s x}}{1-\rme^{\Re\nu^\s T}} \int_0^T \rme^{\Re \nu^\s (T-s)} \,\rmd s\ + \int_0^x \rme^{\Re \nu^\s (x-s)} \,\rmd s\ \right) |P_1 P^\s e_2| \right. \\ & \qquad
+ \left. \left(\frac{\rme^{\Re \nu^\un (x-T)}}{1-\rme^{-\Re\nu^\un T}} \int_T^0 \rme^{-\Re \nu^\un s} \,\rmd s\ + \int_T^x \rme^{\Re \nu^\un (x-s)} \,\rmd s\ \right) |P_1 P^\un e_2| \right] |g|_{X^0} \\
\leq &\ \frac{2}{1-\rme^{-\min(|\Re\nu^\mathrm{s,u}|)T}} \left( \frac{|P_1 P^\s e_2|}{|\Re\nu^\s|} + \frac{|P_1 P^\un e_2|}{|\Re\nu^\un|} \right) |g|_{X^0},
\end{align*}
which establishes (\ref{e:test}).

Finally, for $g\in X^k_\mathrm{per}$, equation (\ref{d:s}) in Lemma~\ref{l:1} shows that $\mathcal{T}$ has the representation given in (\ref{d:t}). Furthermore, proceeding as in the preceding paragraph and using the estimates
\begin{align*}
\|A^{-k}|_{E^\s}\|\leq\frac{1}{|\nu^\s|^k}, \qquad
\|A^{-k}|_{E^\un}\|\leq\frac{1}{|\nu^\un|^k}
\end{align*}
establishes the estimate (\ref{d:test}).
\end{Proof} 

\subsection{Solutions to the linearized eigenvalue problem for the $u$-component}

We fix $T=2\pi/\kappa$ and consider the differential equation
\begin{equation}\label{e:3}
u_{xx} = - \left(\frac{2\rmi}{\alpha}+\omega\right) u_x + \frac{1}{\alpha^2} u + g(x),
\qquad u\in\C, \qquad 0\leq x\leq 2\pi/\kappa
\end{equation}
for $|\alpha|\ll1$. The following result establishes the existence of $2\pi/\kappa$-periodic solutions of (\ref{e:3}). We note that the results of Lemmas~\ref{l:3}-\ref{l:4} and Propositions~\ref{p:1}-\ref{p:2} hold for $|\alpha|\ll 1$. For clarity, we consider the case $\alpha >0$ in the proofs. 
\begin{Lemma}\label{l:3}
Fix $\omega>0$, then there are constants $\alpha_0>0$ and $C_0>0$ so that the following is true. For each $|\alpha|\leq\alpha_0$ and each $g\in X^0$, equation (\ref{e:3}) has a unique $2\pi/\kappa$-periodic solution $u=\mathcal{T}(\alpha)g\in X^2$, and the operator $\mathcal{T}(\alpha)$ satisfies $\|\mathcal{T}(\alpha)\|_{L(X^0,X^0_\mathrm{per})}\leq C_0|\alpha|$. Furthermore, we have
\begin{align}
|\mathcal{T}(\alpha) g|_{X^0_\mathrm{per}} & \leq C_0\alpha^2 |g|_{X^1_\mathrm{per}},
\qquad \mbox{ for } g\in X^1_\mathrm{per} \label{d:talpha1} \\ \label{d:talpha4}
\left|\mathcal{T}(\alpha)g + \alpha^2 g + 2\rmi\alpha^3 g_x\right|_{X^0_\mathrm{per}} & \leq C_0 \alpha^4 |g|_{X^3_\mathrm{per}},
\qquad \mbox{ for } g\in X^3_\mathrm{per}.
\end{align}
\end{Lemma}

\begin{Proof}
First, we rewrite the second order equation~(\ref{e:3}) as the first-order system
\begin{equation} \label{e:3_sys}
U_x = \begin{pmatrix} 0 & 1 \\ \frac{1}{\alpha^2} & -\frac{2i}{\alpha} - \omega \end{pmatrix} U + \begin{pmatrix} 0 \\ 1 \end{pmatrix} g(x)
=: A(\alpha) U + e_2 g(x), \qquad U = \begin{pmatrix} u \\ u_x \end{pmatrix}.
\end{equation}
The eigenvalues $\nu_\pm$ of $A(\alpha)$ are roots of the polynomial $\nu^2 + \left( \frac{2i}{\alpha} + \omega \right) \nu - \frac{1}{\alpha^2}$, which are given by
\begin{align}\label{e:evals}
\Re\nu_\pm = \pm (1+\rmO(\sqrt{\alpha})) \sqrt{\frac{\omega}{2\alpha}}, \qquad
\Im\nu_\pm = \frac{-1+\rmO(\sqrt{\alpha})}{\alpha}.
\end{align}
In particular, there is an $\alpha_0>0$ so that the matrices $A(\alpha)$ are hyperbolic with $|\Re\nu_\pm|\geq\sqrt{\omega/4\alpha}$ for all $0<\alpha\leq\alpha_0$, and Lemma~\ref{l:2} guarantees that (\ref{e:3_sys}) has a unique $2\pi/\kappa$-periodic solution for each $0<\alpha\leq\alpha_0$ and that this solution is given by $u=\mathcal{T}(\alpha)g$.

It remains to verify the estimates for $\mathcal{T}(\alpha)$. For $0<\alpha\ll1$, the stable and unstable eigenspaces of $A(\alpha)$ are given by
\begin{align*}
E^\mathrm{s} = \mathbb{R} \begin{pmatrix} \frac{1}{\nu_-} \\ 1 \end{pmatrix}, \qquad
E^\mathrm{u} = \mathbb{R} \begin{pmatrix} \frac{1}{\nu_+} \\ 1 \end{pmatrix}
\end{align*}
and the associated spectral projections can be written as
\begin{align*}
P^\mathrm{s} U = \left( \frac{1}{\nu_+} - \frac{1}{\nu_-} \right)^{-1} \begin{pmatrix} \frac{1}{\nu_-} \\ 1 \end{pmatrix} 
\left( -1, \frac{1}{\nu_+} \right) U, \qquad
P^\mathrm{u} U = \left( \frac{1}{\nu_-} - \frac{1}{\nu_+} \right)^{-1} \begin{pmatrix} \frac{1}{\nu_+} \\ 1 \end{pmatrix} 
\left( -1, \frac{1}{\nu_-} \right) U.
\end{align*}
We conclude from (\ref{e:evals}) that for all $0<\alpha\ll1$
\begin{align*}
|\nu_\pm| = \frac{1+\rmO(\sqrt{\alpha})}{\alpha}, \quad
\frac{1}{|\nu_\pm|} = \alpha(1+\rmO(\sqrt{\alpha})), \quad
\frac{1}{|\nu_- - \nu_+|} = \frac{\sqrt{\alpha}(1+\rmO(\alpha^2))}{2\sqrt{\omega}}, \quad
\frac{1}{1-\rme^{-2\pi\min(|\Re\nu_\pm|)/\kappa}} \leq 2.
\end{align*}
Hence, we have
\begin{align*}
|P_1 P^\s e_2| & = \left| \frac{1}{\nu_- \nu_-+} \left( \frac{1}{\nu_+} - \frac{1}{\nu_-} \right)^{-1}\right| =
\frac{1}{|\nu_- - \nu_+|} = \frac{\sqrt{\alpha}(1+\rmO(\alpha^2))}{2\sqrt{\omega}}
\end{align*}
and similarly
\begin{align*}
|P_1 P^\un e_2| = \frac{\sqrt{\alpha}(1+\rmO(\alpha^2))}{2\sqrt{\omega}}.
\end{align*}
Substituting these expressions into (\ref{e:test}), we see that there is an $\alpha_0>0$ so that
\begin{align*}
\|\mathcal{T}(\alpha)\|\opno \stackrel{(\ref{e:test})}{\leq} \frac{2}{1-\rme^{-2\pi\min(|\Re\nu_\pm|)/\kappa}} \left( \frac{|P_1 P^\s e_2|}{|\Re\nu^\s|} + \frac{|P_1 P^\un e_2|}{|\Re\nu^\un|} \right)
\leq \frac{4\alpha}{\omega}
\end{align*}
for all $0<\alpha<\alpha_0$ as claimed.

Next, we establish (\ref{d:talpha1}) by using the estimate (\ref{d:test}) in the expression (\ref{d:t}) with $k=1$. Since
\begin{align*}
A^{-1} = \alpha^2 \begin{pmatrix} \frac{2i}{\alpha} + \omega & 1 \\ \frac{1}{\alpha^2} & 0 \end{pmatrix},
\end{align*}
we see that $P_1 A^{-1} e_2=\alpha^2$. Hence, (\ref{d:t})-(\ref{d:test}) with $k=1$ for $g\in X^1_\mathrm{per}$ yield
\begin{align*}
|\mathcal{T}(\alpha) g|_{X^0_\mathrm{per}}
& \stackrel{(\ref{d:t})}{\leq} \left|P_1 A^{-1} \mathcal{S} e_2 g_x\right|_{X^0_\mathrm{per}} + |P_1 A^{-1} e_2| |g| \\
& \stackrel{(\ref{d:test})}{\leq} 4 \left( \frac{|P_1 P^\s e_2|}{|\Re\nu^\s|} \frac{1}{|\nu^s|} + \frac{|P_1 P^\un e_2|}{|\Re\nu^\un|} \frac{1}{|\nu^\un|}\right) \left|g_x\right|_{X^0_\mathrm{per}} + |P_1 A^{-1} e_2| |g|_{X^0_\mathrm{per}} \\
& \leq \frac{8}{\omega} \alpha^2 \left|g_x\right|_{X^0_\mathrm{per}} + \alpha^2 |g|_{X^0_\mathrm{per}}
\end{align*}
as claimed.

Finally, we take $g\in X^3\per$. Calculating $P_1A^{-j}e_2$ iteratively for $j=1,2,3$ shows that
\begin{align*}
P_1 A^{-1} e_2 = \alpha^2, \qquad
P_1 A^{-2} e_2 = 2\rmi\alpha^3 + \omega\alpha^4, \qquad
|P_1 A^{-3} e_2| \leq 5\alpha^4,
\end{align*}
and (\ref{d:t})-(\ref{d:test}) with $k=3$ gives
\begin{align*}
|\mathcal{T}(\alpha) g + \alpha^2 g + 2\rmi\alpha^3 g_x|_{X^0_\mathrm{per}}
& \leq \left|P_1 A^{-3} \mathcal{S} e_2 g_{xxx}\right|_{X^0_\mathrm{per}} + \omega\alpha^4 |g_x|_{X^0_\mathrm{per}} + |P_1 A^{-3} e_2| |g_{xx}|_{X^0_\mathrm{per}} \\
& \leq 4\left( \frac{|P_1 P^\s e_2|}{|\Re\nu^\s|} \frac{1}{|\nu^s|^3} + \frac{|P_1 P^\un e_2|}{|\Re\nu^\un|^3} \frac{1}{|\nu^\un|^3}\right) |g_{xxx}|_{X^0_\mathrm{per}} + (5+\omega) \alpha^4 |g|_{X^2_\mathrm{per}} \\
& \leq C_0 \alpha^4 |g|_{X^3_\mathrm{per}},
\end{align*}
which completes the proof of (\ref{d:talpha4}).
\end{Proof} 

Next, we consider the eigenvalue problem (\ref{evp:1}) for the $u$-component, which is given by the differential equation
\begin{equation}\label{e:4}
u_{xx} = - \left(\frac{2\rmi}{\alpha}+\omega\right) u_x + \frac{1}{\alpha^2} u + (\lambda - f_1(x)) u - f_2(x) v
\qquad u\in\C, \qquad 0\leq x\leq 2\pi/\kappa.
\end{equation}
We will vary $\lambda$ in the rectangle $\Lambda_R$ of width $2R$ and height one that we defined in (\ref{d:1}). For fixed functions $f_1,f_2\in X^0_\mathrm{per}$, we define the multiplication operators $B_1(\lambda), B_2\in L(X^0)$ via $B_1(\lambda)v:=(\lambda - f_1(x))v$ and $B_2v:=-f_2v$ and note that the map
\begin{align*}
B_1: \Lambda_R\longrightarrow L(X^0), \quad \lambda\longmapsto B_1(\lambda)
\end{align*}
is analytic. The next result establishes the existence of $2\pi/\kappa$-periodic solutions of (\ref{e:4}).

\begin{Proposition}\label{p:1}
For each fixed choice of $R>0$, $\omega>0$, and $f_1,f_2\in X^0_\mathrm{per}$, there are constants $\alpha_1>0$ and $C_1>0$ so that the following is true. For each $(\lambda,\alpha)\in\Lambda_R\times[-\alpha_1,\alpha_1]\setminus\{0\}$ and $v\in X^0$, equation (\ref{e:4}) has a unique $2\pi/\kappa$-periodic solution $u=\mathcal{D}(\lambda,\alpha)g\in X^2$, where 
\begin{align*}
\mathcal{D}(\lambda,\alpha)=(1-\mathcal{T}(\alpha)B_1(\lambda))^{-1}\mathcal{T}(\alpha)B_2.
\end{align*}
The function
\begin{align*}
\mathcal{D}(\cdot,\alpha):\quad \Lambda_R\longrightarrow L(X^0,X^0_\mathrm{per}), \quad \lambda\longmapsto \mathcal{D}(\lambda,\alpha)
\end{align*}
is well defined and analytic in $\lambda$ for each fixed $\alpha \in [-\alpha_0,\alpha_0]\setminus\{0\}$, and we have
\begin{align}\label{e:d}
\|\mathcal{D}(\lambda,\alpha)\|_{L(X^0,X^0_\mathrm{per})} \leq C_1\alpha, \qquad
\|\mathcal{D}_\lambda(\lambda,\alpha)\|_{L(X^0,X^0_\mathrm{per})} \leq C_1\alpha^2.
\end{align}
Finally, setting $\mathcal{D}(\lambda,0):=0$, the mapping $\mathcal{D}:\Lambda_R\times[-\alpha_1,\alpha_1]\to L(X^0,X^0_\mathrm{per})$ is analytic in $\lambda$, and $\mathcal{D}$ and its derivatives in $\lambda$ are continuous in $\alpha$.
\end{Proposition}

\begin{Proof}
The eigenvalue problem (\ref{e:4}) coincides with (\ref{e:3}) provided we set $g(x)=B_1(\lambda)u+B_2v$ in (\ref{e:3}). Lemma~\ref{l:3} therefore shows that there is an $\alpha_0>0$ so that (\ref{e:4}) has a solution $u\in X^2$ for some $\alpha\in(0,\alpha_0]$ if and only $u\in X^0_\mathrm{per}$ satisfies
\begin{equation}\label{e:fpu}
u = \mathcal{T}(\alpha) \left( B_1(\lambda)u + B_2 v \right),
\end{equation}
where $\mathcal{T}(\alpha)$ has been defined in Lemma~\ref{l:3}. We set $M:=\max(|f_1|_{X^0_\mathrm{per}}+R,|f_2|_{X^0_\mathrm{per}})$ and note that Lemma~\ref{l:3} implies that there is a constant $C_0$ so that the estimates
\begin{align*}
\|\mathcal{T}(\alpha) B_1(\lambda)\|\xper \leq C_0 M \alpha, \qquad
\|\mathcal{T}(\alpha) B_2\|\xper \leq C_0 M \alpha
\end{align*}
hold uniformly in $(\lambda,\alpha)\in\Lambda_R\times(0,\alpha_0]$. Choosing $\alpha_1:=\min(\alpha_0,\frac{1}{2C_0M})$, we see that $\|\mathcal{T}(\alpha)B_1(\lambda)\|\xper\leq\frac12$ for $(\lambda,\alpha)\in\Lambda_R\times(0,\alpha_1]$, and the fixed-point equation (\ref{e:fpu}) therefore has a unique solution
\begin{align*}
u = (1 - \mathcal{T}(\alpha) B_1(\lambda))^{-1} \mathcal{T}(\alpha) B_2 v =: \mathcal{D}(\lambda,\alpha) v
\end{align*}
in $X^0_\mathrm{per}$ for each $(\lambda,\alpha)\in\Lambda_R\times(0,\alpha_1]$ and $v\in X^0$. Furthermore, we have
\begin{align*}
\mathcal{D}(\lambda,\alpha) = (1 - \mathcal{T}(\alpha) B_1(\lambda))^{-1} \mathcal{T}(\alpha) B_2 \in L(X^0,X^0_\mathrm{per})
\end{align*}
with
\begin{align*}
\|\mathcal{D}(\lambda,\alpha)\|\opno = \|(1-\mathcal{T}(\alpha)B_1(\lambda))^{-1}\|\opno \|\mathcal{T}(\alpha)B_2\|\xper \leq 2 \|\mathcal{T}(\alpha)B_2\|\xper \leq 2C_0M \alpha
\end{align*}
for $(\lambda,\alpha)\in\Lambda_R\times(0,\alpha_1]$.

Analyticity of $\mathcal{D}(\lambda,\alpha)$ in $\lambda$ follows from analyticity of $B_1(\lambda)$. To estimate the derivative $\mathcal{D}_\lambda(\lambda,\alpha)$, we differentiate the identity
\begin{align*}
(1 - \mathcal{T}(\alpha) B_1(\lambda)) \mathcal{D}(\lambda,\alpha) = \mathcal{T}(\alpha) B_2
\end{align*}
with respect to $\lambda$ to get
\begin{align*}
- \mathcal{T}(\alpha) \frac{\rmd B_1}{\rmd\lambda}(\lambda)) \mathcal{D}(\lambda,\alpha)
+ (1 - \mathcal{T}(\alpha) B_1(\lambda)) \mathcal{D}_\lambda(\lambda,\alpha) = 0.
\end{align*}
Since $\frac{\rmd B_1}{\rmd\lambda}(\lambda)=1$, we find
\begin{align*}
\mathcal{D}_\lambda(\lambda,\alpha) = (1 - \mathcal{T}(\alpha) B_1(\lambda))^{-1} \mathcal{T}(\alpha) \mathcal{D}(\lambda,\alpha)
\end{align*}
and therefore
\begin{align*}
\|\mathcal{D}_\lambda(\lambda,\alpha)\|\opno \leq
\|\left( 1 - \mathcal{T}(\alpha)B_1(\lambda) \right)^{-1}\|\opno \|\mathcal{T}(\alpha)\|\xper \|\mathcal{D}(\lambda,\alpha)\|\opno \leq 4C_0^2 M \alpha^2
\end{align*}
as claimed. Note that $\mathcal{D}(\lambda,\alpha)$ and its derivatives in $\lambda$ converge to zero in $L(X^0,X^0_\mathrm{per})$ as $\alpha\searrow0$, and we can therefore extend $\mathcal{D}(\lambda,\alpha)$ into $\alpha=0$ by setting $\mathcal{D}(\lambda,0)=0$.
\end{Proof} 

\subsection{Solutions to the linearized eigenvalue problem for the $v$-component for $\delta=0$}

It remains to find nontrivial $2\pi/\kappa$-periodic solutions of the eigenvalue problem (\ref{evp:2}), given by
\begin{equation}\label{e:5}
\lambda v = \omega v_x + \bar{g} v + g_1(x) u
\end{equation}
where $u=\mathcal{D}(\lambda,\alpha)v$ is the unique $2\pi/\kappa$-periodic solution of (\ref{e:4}) for a given function $v$, whose existence was shown in Proposition~\ref{p:1}. Our first result shows existence of solutions to (\ref{e:5}) without enforcing periodicity in $x$.  

\begin{Lemma}\label{l:4}
Fix $R>0$, $\omega>0$, $\bar{g}\in\R$, and $f_1,f_2,g_1\in X^0_\mathrm{per}$, then there are constants $\alpha_2>0$ and $C_2>0$ so that the following is true. For each $(\lambda,\alpha)\in\Lambda_R\times [-\alpha_2,\alpha_2]\setminus\{0\}$, equation (\ref{e:5}) with $u=\mathcal{D}(\lambda,\alpha)v$ has a nontrivial solution $v=\mathcal{C}(\lambda,\alpha)\in X^1$, and this solution is unique up to a constant factor. In addition:
\begin{compactitem}
\item We have
\begin{align}
\mathcal{C}(\lambda,\alpha) & = (1+\mathcal{Q}(\lambda,\alpha))^{-1} \rme^{(\lambda-\bar{g})x/\omega}
\label{def:c} \\ \label{def:q}
[\mathcal{Q}(\lambda,\alpha) h](x) & := \frac{1}{\omega} \int_0^x \rme^{(\lambda-\bar{g})(x-y)/\omega} g_1(y) [\mathcal{D}(\lambda,\alpha)h](y)\, \rmd y, \qquad h\in X^0.
\end{align}
\item For each fixed $\alpha\in[-\alpha_2,\alpha_2]\setminus\{0\}$, the mappings
\begin{align*}
\mathcal{Q}(\cdot,\alpha):\; \Lambda_R\to L(X^0),\; \lambda\mapsto \mathcal{Q}(\lambda,\alpha), \qquad
\mathcal{C}(\cdot,\alpha):\; \Lambda_R\to X^0,\; \lambda\mapsto \mathcal{C}(\lambda,\alpha)
\end{align*}
are analytic in $\lambda$ with
\begin{align}\label{e:cq}
\left\|\frac{\rmd^k}{\rmd\lambda^k}\mathcal{Q}(\lambda,\alpha)\right\|_{L(X^0)} \leq C_2\alpha, \quad
\left|\frac{\rmd^k}{\rmd\lambda^k}\mathcal{C}(\lambda,\alpha)\right|_{X^0} \leq C_2, \quad
\left|\mathcal{C}(\lambda,\alpha) - \rme^{(\lambda-\bar{g})x/\omega}\right|_{X^0} \leq C_2\alpha
\end{align}
for $k=0,1$ uniformly in $(\lambda,\alpha)\in\Lambda_R\times(0,\alpha_2]$. Furthermore, 
setting $\mathcal{Q}(\lambda,0)=0$ and $\mathcal{C}(\lambda,0)=0$ shows that $\mathcal{Q}$ and $\mathcal{C}$ are analytic in $\lambda$ and continuous in $\alpha$ on $\Lambda_R\times[-\alpha_2,\alpha_2]$.
\end{compactitem}
\end{Lemma}

\begin{Proof}
Let $\alpha_1$ be the constant from Proposition~\ref{p:1}. For each $0<\alpha\leq\alpha_1$, a function $v\in X^1$ is then a solution of (\ref{e:5}) for $u=\mathcal{D}(\lambda,\alpha)v$ with $v(0)=1$ if and only if $v\in X^0$ satisfies the fixed-point equation
\begin{align*}
v(x) = \rme^{(\lambda-\bar{g})x/\omega} - \frac{1}{\omega} \int_0^x \rme^{(\lambda-\bar{g})(x-y)/\omega}g_1(y)[\mathcal{D}(\lambda,\alpha)v](y)\,\rmd y =: \rme^{(\lambda-\bar{g})x/\omega} - [\mathcal{Q}(\lambda,\alpha) v](x).
\end{align*}
Using the estimates for $\mathcal{D}(\lambda,\alpha)$ from Proposition~\ref{p:1}, we have
\begin{align*}
\|\mathcal{Q}(\lambda,\alpha)\|_{L(X^0)}
\leq \frac{1}{\omega} \int_0^{2\pi/\kappa} |\rme^{(\lambda-\bar{g})(2\pi/\kappa - y)/\omega}|\,\rmd y\ \|\mathcal{D}(\lambda,\alpha)\|_{L(X^0)}\ |g_1|\xper
\leq \frac{2\pi}{\omega\kappa} \rme^{2\pi(R+|\bar{g}|)/(\omega\kappa)} C_1 \alpha
=: \tilde{C}_1 \alpha.
\end{align*}
In particular, setting $\alpha_2:=\min(\alpha_1,\frac{1}{2\tilde{C}_1})$, we can solve the fixed-point equation
\begin{align*}
v = \rme^{(\lambda-\bar{g})x/\omega} - \mathcal{Q}(\lambda,\alpha) v
\end{align*}
uniquely for $v$ for each $(\lambda,\alpha)\in\Lambda_R\times(0,\alpha_2]$ to get
\begin{align*}
v = \mathcal{C}(\lambda,\alpha) := \left(1+\mathcal{Q}(\lambda,\alpha)\right)^{-1} \rme^{(\lambda-\bar{g})x/\omega}.
\end{align*}
Furthermore,
\begin{align*}
|\mathcal{C}(\lambda,\alpha)|_{X^0} \leq 2\rme^{2\pi(R+|\bar{g}|)/(\kappa\omega)} =: \tilde{C}_2, \qquad
|\mathcal{C}(\lambda,\alpha) - \rme^{(\lambda-\bar{g})x/\omega}|_{X^0} = |\mathcal{Q}(\lambda,\alpha) \mathcal{C}(\lambda,\alpha)|_{X^0} \leq \tilde{C}_1 \tilde{C}_2 \alpha,
\end{align*}
and setting $C_2:=\max(\tilde{C}_1,\tilde{C}_2,\tilde{C}_1\tilde{C}_2)$ completes the proof of the estimates. Finally, analyticity of $\mathcal{D}(\lambda,\alpha)$ implies analyticity of $\mathcal{Q}(\lambda,\alpha)$ and $\mathcal{C}(\lambda,\alpha)$, and the estimates for the derivatives in $\lambda$ follow in a similar fashion.
\end{Proof}

Our next result characterizes the set of $\lambda\in\Lambda_R$ for which the eigenvalue problem
\begin{align}
\lambda u & = u_{xx} + \left(\frac{2\rmi}{\alpha}+\omega\right) u_x - \frac{1}{\alpha^2} u + f_1(x) u + f_2(x) v
\label{evp:5} \\ \label{evp:6}
\lambda v & = \omega v_x + \bar{g} v + g_1(x) u
\end{align}
has a nontrivial $2\pi/\kappa$-periodic solution for an appropriate $|\alpha|\ll1$.

\begin{Proposition}\label{p:2}
Fix $\omega>0$, $\bar{g}\in\R$, $R>|\bar{g}|$, and $f_1,f_2,g_1\in X^3_\mathrm{per}$, then there are constants $\alpha_3>0$ and $C_3>0$ so that the following is true. There is a unique function $\lambda_0^*:[-\alpha_3,\alpha_3]\to\Lambda_R$ so that (\ref{evp:5})-(\ref{evp:6}) has a nontrivial $2\pi/\kappa$-periodic solution for $(\lambda,\alpha)\in\Lambda_R\times[-\alpha_3,\alpha_3]$ if and only if $\lambda=\lambda_0^*(\alpha)$. Furthermore, upon defining the real constants
\begin{align*}
\lambda_2 := \frac{\kappa}{2\pi} \int_0^{2\pi/\kappa} f_2(x) g_1(x)\,\rmd x, \qquad
\lambda_3 := \frac{\kappa}{2\pi} \int_0^{2\pi/\kappa} f_2^\prime(x) g_1(x)\,\rmd x,
\end{align*}
we have the expansion
\begin{align*}
|\lambda_0^*(\alpha) - (\bar{g} + \lambda_2\alpha^2 + 2\rmi\lambda_3\alpha^3)| \leq C_3\alpha^4
\end{align*}
for the dispersion curve, and the associated eigenfunctions $(u,v)(x;\alpha)$ satisfy
\begin{align*}
\left|u(x;\alpha) - (\alpha^2 f_2(x) + 2\rmi \alpha^3 f_2^\prime(x))\right|_{X^0_\mathrm{per}} \leq C_3\alpha^4, \qquad
\left|v(x;\alpha) - (1 + \alpha^2 v_2(x) + 2\rmi \alpha^3 v_3(x))\right|_{X^0_\mathrm{per}} \leq C_3\alpha^4,
\end{align*}
where
\begin{align*}
v_2(x) = \frac{1}{\omega} \left(\lambda_2 x - \int_0^x f_2(y) g_1(y)\, \rmd y \right), \qquad
v_3(x) = \frac{1}{\omega} \left(\lambda_3 x - \int_0^x f_2^\prime(y) g_1(y)\, \rmd y \right).
\end{align*}
\end{Proposition}

\begin{Proof}
We will use the notation and constants introduced in Lemma~\ref{l:4}. Proposition~\ref{p:2} and Lemma~\ref{l:4} show that (\ref{evp:5})-(\ref{evp:6}) has a nontrivial solution $(u,v)\in X^2_\mathrm{per}\times X^1_\mathrm{per}$ for some $(\lambda,\alpha)\in\Lambda_R\times(0,\alpha_2]$ if and only if $(u,v)=(\mathcal{D}(\lambda,\alpha)\mathcal{C}(\lambda,\alpha),\mathcal{C}(\lambda,\alpha))\in X^0_\mathrm{per}\times X^0_\mathrm{per}$ with $v(2\pi/\kappa)=v(0)$. It therefore suffices to solve $\Delta(\lambda,\alpha):=v(2\pi/\kappa)-v(0)=0$. Using the expression for $\mathcal{C}(\lambda,\alpha)$ from Lemma~\ref{l:4}, we can write $\Delta(\lambda,\alpha)$ as
\begin{align}
\Delta(\lambda,\alpha)
& = \rme^{2\pi(\lambda-\bar{g})/(\kappa\omega)} - 1 - \frac{1}{\omega} \int_0^{2\pi/\kappa} \rme^{(\lambda-\bar{g})(2\pi/\kappa-y)/\omega} g_1(y) [\mathcal{D}(\lambda,\alpha)v](y)\,\rmd y \nonumber \\ \label{e:Delta0}
& = \rme^{2\pi(\lambda-\bar{g})/(\kappa\omega)} - 1 - [\mathcal{Q}(\lambda,\alpha) \mathcal{C}(\lambda,\alpha)]\left(2\pi/\kappa\right)
\end{align}
and conclude that $\Delta(\lambda,\alpha)$ is analytic in $\lambda$. Lemma~\ref{l:4} shows that
\begin{align*}
\left|[\mathcal{Q}(\alpha,\lambda) \mathcal{C}(\alpha,\lambda)]\left(2\pi/\kappa\right)\right| \leq \|\mathcal{Q}(\alpha,\lambda)\|_{L(X^0)} |\mathcal{C}(\alpha,\lambda)|_{X^0} \leq C_2^2 \alpha.
\end{align*}
Using the extensions of $\mathcal{Q}$ and $\mathcal{C}$ into $\alpha=0$, it follows that
\begin{align*}
\Delta(\lambda,0) = \rme^{2\pi(\lambda-\bar{g})/(\kappa\omega)} - 1
\end{align*}
so that $\lambda=\bar{g}\in\Lambda_R$ is the only solution in $\Lambda_R$ when $\alpha=0$. Since $\Delta_\lambda(\bar{g},0)=2\pi/(\kappa\omega)\neq0$, we can solve $\Delta(\lambda,\alpha)=0$ near $(\bar{g},0)$ uniquely by the implicit function theorem and conclude that there is an $\alpha_3>0$ and a unique function $\lambda^0_*:[0,\alpha_3]\to\C$ with $\lambda^0_*(0)=\bar{g}$ so that $\Delta(\lambda,\alpha)=0$ for $(\lambda,\alpha)\in\Lambda_R\times[0,\alpha_3]$ if and only if $\lambda=\lambda^0_*(\alpha)$.

It remains to establish the expansions for $\lambda=\lambda^0_*(\alpha)$ and the associated solutions $(u,v):=(u,v)(x;\alpha)$ of (\ref{evp:5})-(\ref{evp:6}). Setting $B_1(\lambda)v:=(\lambda-f_1(x))v$ and $B_2v:=-f_2v$, we can represent $(u,v,\lambda)$ via
\begin{align}
u & = \mathcal{D}(\lambda,\alpha)v = (1-\mathcal{T}(\alpha)B_1(\lambda))^{-1}\mathcal{T}(\alpha) B_2 v
\label{r:u} \\ \label{r:v}
v & = \rme^{(\lambda-\bar{g})x/\omega} - \frac{1}{\omega} \int_0^x \rme^{(\lambda-\bar{g})(x-y)/\omega} g_1(y) u(y)\,\rmd y 
\\ \label{r:l}
0 & = \rme^{2\pi(\lambda-\bar{g})/(\kappa\omega)} - 1 - \frac{1}{\omega} \int_0^{2\pi/\kappa} \rme^{(\lambda-\bar{g})(2\pi/\kappa-y)/\omega} g_1(y) u(y)\,\rmd y
\end{align}
where $\lambda$ is always evaluated at $\lambda^0_*(\alpha)$. Throughout, we will denote by $C$ a constant that depends only on $\omega$, $\bar{g}$, $R$, $|f_1|_{X^3_\mathrm{per}}$, $|f_2|_{X^3_\mathrm{per}}$, and $|g_1|_{X^3_\mathrm{per}}$ and that may change from estimate to estimate. We will use Landau symbols only when we can estimate them by such a constant multiplied by the argument of the Landau symbol.

First, equation (\ref{e:cq}) shows that $v=\mathcal{C}(\lambda_*(\alpha),\alpha)=\rmO(1)$. Equation (\ref{r:u}) and the estimate (\ref{e:d}) for $\mathcal{D}$ then imply that $u=\rmO(\alpha)$, and using this estimate in (\ref{r:l}) proves that $\lambda^0_*(\alpha)=\bar{g}+\rmO(\alpha)$. Using these expansions in (\ref{r:v}), we conclude that
\begin{align*}
v = 1 + \rmO(\alpha), \qquad
v_x = \frac{\lambda-\bar{g}}{\omega} v - \frac{g_1(x)}{\omega} u = \rmO(\alpha),
\end{align*}
and we therefore have $|v|_{X^1_\mathrm{per}}\leq C$ for a constant $C$ as above. Using the estimate $|\mathcal{T}(\alpha) g|_{X^0_\mathrm{per}} \leq C_0\alpha^2 |g|_{X^1_\mathrm{per}}$ from Lemma~\ref{l:3}, we obtain from (\ref{r:u}) that
\begin{align*}
|u|\xper
\leq \|(1-\mathcal{T}(\alpha)B_1(\lambda))^{-1}\|\opno |\mathcal{T}(\alpha)B_2v|\xper
\leq 2 C_0 \alpha^2 |B_2 v|_{X^1_\mathrm{per}}
\leq 2 C_0 C \alpha^2,
\end{align*}
which shows that $u=\rmO(\alpha^2)$. Using this estimate for $u$ in (\ref{r:l}) gives
\begin{align}\label{e:expl0}
\lambda^0_*(\alpha) =: \bar{g} + \lambda_2(\alpha), \qquad \lambda_2(\alpha)=\rmO(\alpha^2)
\end{align}
and therefore
\begin{align*}
v_x = \frac{\lambda-\bar{g}}{\omega} v - \frac{g_1(x)}{\omega} u = \rmO(\alpha^2).
\end{align*}
In particular, we can write $v=1+\tilde{v}$ with $|\tilde{v}|_{X^1_\mathrm{per}}\leq C\alpha^2$.

Next, we verify the expansion for $u$. We write
\begin{align}\label{e:expu0}
u = \alpha^2 f_2 + 2\rmi \alpha^3 f_2^\prime + \tilde{u}
\end{align}
and need to show that $\tilde{u}=\rmO(\alpha^4)$. Equation (\ref{r:u}) can be written as $u=\mathcal{T}(\alpha)B_1(\lambda)u-\mathcal{T}(\alpha)(f_2v)$ and substituting our expression for $u$ as well as $v=1+\tilde{v}$ into this equation gives
\begin{align}\label{e:expu1}
\alpha^2 f_2 + 2\rmi \alpha^3 f_2^\prime + \tilde{u} = 
\alpha^2 \mathcal{T}(\alpha)B_1(\lambda) f_2
+ 2\rmi \alpha^3 \mathcal{T}(\alpha)B_1(\lambda) f_2^\prime
+ \mathcal{T}(\alpha)B_1(\lambda) \tilde{u}
- \mathcal{T}(\alpha)f_2 - \mathcal{T}(\alpha)(f_2\tilde{v}).
\end{align}
Since $B_1(\lambda)f_2$ and $B_1(\lambda)f_2^\prime$ are in $X^1_\mathrm{per}$ and $|f_2\tilde{v}|_{X^1_\mathrm{per}}=\rmO(\alpha^2)$, we can use the estimate (\ref{d:talpha1}) to get
\begin{align*}
\left|\alpha^2 \mathcal{T}(\alpha)B_1(\lambda) f_2
+ 2\rmi \alpha^3 \mathcal{T}(\alpha)B_1(\lambda) f_2^\prime
+ \mathcal{T}(\alpha)(f_2\tilde{v})\right|_{X^0_\mathrm{per}} \leq C\alpha^4.
\end{align*}
Furthermore, since $f_2\in X^3_\mathrm{per}$, equation (\ref{d:talpha4}) shows that
\begin{align*}
\left|\mathcal{T}(\alpha)f_2 + \alpha^2 f_2 + 2\rmi \alpha^3 f_2^\prime\right|_{X^0_\mathrm{per}} \leq C\alpha^4.
\end{align*}
Hence, (\ref{e:expu1}) becomes
\begin{align*}
\tilde{u} = \mathcal{T}(\alpha)B_1(\lambda) \tilde{u} + \rmO(\alpha^4)
\end{align*}
in $X^0_\mathrm{per}$ and since $\|\mathcal{T}(\alpha)B_1(\lambda)\|_{L(X^0)}\leq\frac12$, we conclude that $|\tilde{u}|_{X^0}\leq C\alpha^4$ as claimed.

Our last step is to verify the expansions for $v$ and $\lambda^0_*(\alpha)$. We substitute the expansions (\ref{e:expl0}) and (\ref{e:expu0}) into (\ref{r:v}) to get
\begin{align*}
v(x)
& = \rme^{(\lambda-\bar{g})x/\omega} - \frac{1}{\omega} \int_0^x \rme^{(\lambda-\bar{g})(x-y)/\omega} g_1(y) u(y)\,\rmd y \\
& = \rme^{\lambda_2(\alpha)x/\omega} - \frac{1}{\omega} \int_0^x \rme^{\lambda_2(\alpha)(x-y)/\omega} g_1(y)
\left(\alpha^2 f_2(y) + 2\rmi \alpha^3 f_2^\prime(y) + \rmO(\alpha^4)\right)\,\rmd y \\
& = \left[1+\frac{\lambda_2(\alpha)}{\omega}x \right] - \frac{\alpha^2}{\omega} \int_0^x g_1(y) f_2(y) \,\rmd y - \frac{2\rmi\alpha^3}{\omega} \int_0^x g_1(y) f_2^\prime(y) \,\rmd y\ + \rmO(\alpha^4).
\end{align*}
We know that $v(0)=v(2\pi/\kappa)$ which shows that
\begin{align*}
\lambda_2(\alpha)
= \frac{\alpha^2 \kappa}{2\pi} \int_0^{2\pi/\kappa} g_1(y) f_2(y)\,\rmd y\ + \frac{2\rmi\alpha^3 \kappa}{2\pi} \int_0^{2\pi/\kappa} g_1(y) f_2^\prime(y) \, \rmd y\ + \rmO(\alpha^4)
=: \lambda_2\alpha^2 + 2\rmi\lambda_3\alpha^3 + \rmO(\alpha^4)
\end{align*}
as claimed. Substituting this expression into the expansion for $v$ gives
\begin{align*}
v(x)
& = 1 + \frac{\alpha^2}{\omega} \left( \lambda_2 x - \int_0^x g_1(y) f_2(y) \,\rmd y\ \right) + \frac{2\rmi\alpha^3}{\omega} \left( \lambda_3 x - \int_0^x g_1(y) f_2^\prime(y) \,\rmd y \right) + \rmO(\alpha^4).
\end{align*}
This completes the proof of Proposition~\ref{p:2} and the construction of the spectral expansion for the case $\delta=0$.
\end{Proof}

We emphasize that the results in this section extend immediately to the case where the coefficient $\bar{g}$ is a $2\pi/\kappa$-periodic function and not a constant. 

\subsection{Periodic solutions of linear second-order equations with center directions}

In preparation for the case of $\delta >0$, we seek to find nontrivial periodic solutions to second-order equations of the form
\begin{equation} \label{e:genB}
v_{xx} = b_{21}v + b_{22} v_x + h(x), \qquad v \in \mathbb{C}, \qquad 0 \leq x \leq T,
\end{equation}
where the associated first-order system is not guaranteed to be hyperbolic.

\begin{Lemma} \label{l:7}
Fix $T>0$ and coefficients $b_{21},b_{22}$. Assume that there is an $r>0$ such that the eigenvalues $\eta^\mathrm{c},\eta^\s$ of the matrix $B := \begin{pmatrix} 0 & 1 \\ b_{21} & b_{22} \end{pmatrix}$ satisfy $\Re\eta^\s<-2r$ and $|\Re\eta^\mathrm{c}|\leq r$, then the following is true. For each $h \in X^0$, equation~(\ref{e:genB}) has a $T$-periodic solution $v(x)\in X^2_\mathrm{per}$ if and only if (i) there exists an $a^\mathrm{c}\in E^\mathrm{c}$ so that
\begin{align}\label{e:1d}
v(x) = P_1 \left(\rme^{\eta^\mathrm{c}x} a^{\mathrm{c}} + [J h](x)\right), \qquad 0\leq x\leq T
\end{align}
where
\begin{align*}
[J h](x) & = \int_0^x \rme^{\eta^\mathrm{c}(x-s)} P^\mathrm{c} e_2 h(s) \,\rmd s 
+ \frac{\rme^{\eta^\mathrm{s}x}}{1-\rme^{\eta^\mathrm{s}T}} \int_0^T \rme^{\eta^\mathrm{s}(T-s)}P^\s e_2 h(s) \,\rmd s\ + \int_0^x \rme^{\eta^\mathrm{s}(x-s)}P^\s e_2 h(s) \, \rmd s
\end{align*}
and (ii)
\begin{align}\label{e:1D}
(\rme^{\eta^\mathrm{c}T}-1) a^\mathrm{c} + [P^\mathrm{c}Jh](2\pi) = 0.
\end{align}
\end{Lemma}

\begin{Proof}
We write (\ref{e:genB}) as the system
\begin{equation} \label{e:genBsys}
V_x = \begin{pmatrix} 0 & 1 \\ b_{21} & b_{22} \end{pmatrix} V + \begin{pmatrix} 0 \\ h(x) \end{pmatrix} := B V + e_2 h(x), \qquad V = \begin{pmatrix} v \\ v_x \end{pmatrix}.
\end{equation}
Enforcing periodicity in the stable direction as in Lemma~\ref{l:1} and using the variation-of-constants formula in the center direction establishes the expression (\ref{e:1d}) of the solution $v(x)$ in $X^2$. The solution $v$ lies in $X^2_\mathrm{per}$ if and only if $V(T)-V(0)=0$. Since $P^\s(V(T)-V(0))=0$ by construction, this condition reduces to $P^\mathrm{c}(V(T)-V(0))=0$, and substituting the expressions for $V(0)$ and $V(T)$ gives (\ref{e:1D}).
\end{Proof}

\subsection{Solutions to the linearized eigenvalue problem for $\delta >0$}

The next result focuses on the eigenvalue problem
\begin{align}
\lambda u &= u_{xx} + \left( \frac{2 \rmi}{\alpha} + \omega \right) u_x - \frac{1}{\alpha^2} u + f_1(x) u + f_2(x) v
\label{evp:7} \\ \label{evp:8}
\lambda v &= \delta v_{xx} + \left( \frac{2 \rmi \delta}{\alpha} + \omega \right) v_x - \frac{\delta}{\alpha^2} v + g_1(x) u + \bar{g} v
\end{align}
for $0<\delta\ll1$. For each fixed $\delta>0$, we will identify values of $(\lambda,\alpha)$ for which (\ref{evp:7})-(\ref{evp:8}) has a nontrivial $2\pi/\kappa$-periodic solution.

\textbf{Theorem 2.}
\emph{Fix $\omega>0$ and $\bar{g}\in\R$, and assume that $f_1,f_2,g_1$ are given $2\pi/\kappa$-periodic functions of class $C^4$, then there are constants $\delta_0,s_0>0$ and functions $\alpha^*,\lambda^*:[-s_0,s_0]\times(0,\delta_0]\to\C$ with
\begin{align}\label{e:lambda_sd_eqn}
\alpha^*(s,\delta) & = \frac12 \left( s+ \sqrt{s^2+4\sqrt{\delta}} \right) \\
\lambda^*(s,\delta) & = \left\{ \begin{array}{lcl} 
\lambda_*(|s|) + \rmO(\sqrt{\delta}) & \quad & s\geq 0 \\
\bar{g} - \frac{s^2}{\omega} (1+\rmO(\delta^{\frac14})) + \rmO(\sqrt{\delta}) & \quad & s\leq 0
\end{array} \right.
\end{align}
so that (\ref{evp:7})-(\ref{evp:8}) has a nontrivial $2\pi/\kappa$-periodic solution when $(\alpha,\lambda)=(\alpha^*,\lambda^*)(s,\delta)$.
}

\begin{Proof}
Fix $R\geq|\bar{g}|$ and let $\lambda$ vary in $\Lambda_R$. Throughout, we will denote by $C$ a constant that depends only on $\omega$, $\bar{g}$, $R$, $|f_1|_{X^3_\mathrm{per}}$, $|f_2|_{X^3_\mathrm{per}}$, and $|g_1|_{X^3_\mathrm{per}}$ and that may change from estimate to estimate. We will use Landau symbols only when we can estimate them by such a constant multiplied by the argument of the Landau symbol. 

We showed in Proposition~\ref{p:1} that for each given $v\in X^0$ the unique solution of (\ref{evp:7}) is given by $u=\mathcal{D}(\lambda,\alpha)v\in X^2_\mathrm{per}$, and it therefore suffices to solve (\ref{evp:8}) with $u=\mathcal{D}(\lambda,\alpha)v$. Throughout the proof, we will use the scaling $\delta=\alpha^2\beta^2$ with $0<\alpha,\beta\leq1$. These results hold for $|\alpha|\ll 1$, and as before, the proofs consider the case $\alpha >0$ for clarity. With these definitions, (\ref{evp:8}) becomes
\begin{equation} \label{evp:9}
v_{xx} = -\left( \frac{2\rmi\alpha\beta^2+\omega}{\alpha^2\beta^2} \right) v_x + \left(\frac{\lambda+\beta^2-\bar{g}}{\alpha^2 \beta^2} \right) v - \frac{1}{\alpha^2\beta^2} g_1(x) [\mathcal{D}(\lambda,\alpha)v](x).
\end{equation}

We will use Lemma~\ref{l:7} to reformulate (\ref{evp:9}) as a fixed-point problem. The roots of the characteristic equation
\begin{align*}
\eta^2 = -\left( \frac{2\rmi\alpha\beta^2+\omega}{\alpha^2\beta^2} \right) \eta + \left(\frac{\lambda+\beta^2-\bar{g}}{\alpha^2 \beta^2} \right)
\end{align*}
associated with (\ref{evp:9}) are given by
\begin{align}\label{e:nu}
\eta^\mathrm{c}(\lambda,\alpha,\beta) = \frac{\lambda+\beta^2-\bar{g}}{\omega} + \rmO(\alpha\beta^2), \qquad
\eta^\mathrm{s}(\lambda,\alpha,\beta) = -\frac{\omega+\rmO(\alpha\beta^2)}{\alpha^2\beta^2}.
\end{align}
These expressions are analytic in $\lambda\in\Lambda_R$ with $|\Re\eta^\mathrm{c}|\leq C$, $\Re\eta^\mathrm{s}\leq\frac{-\omega}{2\alpha^2\beta^2}$, and $|\eta^\mathrm{c,s}_\lambda(\lambda,\alpha,\beta)|\leq C$ uniformly in $\Lambda_R$. Hence, there is an $\alpha_4$ so that the assumptions of Lemma~\ref{l:7} are met for $0<\alpha\leq\alpha_4$. We conclude that (\ref{evp:9}) has a nontrivial solution $v\in X^2$ if and only if
\begin{align}\label{e:fpe1}
v = \rme^{\eta^\mathrm{c}(\lambda,\alpha,\beta)x} P_1 a^\mathrm{c} - \mathcal{J}(\lambda,\alpha,\beta) B_3 \mathcal{D}(\lambda,\alpha) v, \qquad
\mathcal{J}(\lambda,\alpha,\beta) := \frac{1}{\alpha^2\beta^2} P_1 J(\lambda,\alpha,\beta)
\end{align}
for some $a^\mathrm{c}\in E^\mathrm{c}$, where $B_3\in L(X^0)$ is given by $[B_3v](x)=g_1(x)v(x)$, the operator $J(\lambda,\alpha,\beta)\in L(X^0)$ is defined by
\begin{align*}
[J(\lambda,&\alpha,\beta) h](x)
= \int_0^x \rme^{\eta^\mathrm{c}(\lambda,\alpha,\beta)(x-s)} P^\mathrm{c}(\lambda,\alpha,\beta) e_2 h(s) \,\rmd s \\ &
+ \frac{\rme^{\eta^\mathrm{s}(\lambda,\alpha,\beta)x}}{1-\rme^{2\pi\eta^\mathrm{s}(\lambda,\alpha,\beta)/\kappa}} \int_0^{2\pi/\kappa} \rme^{\eta^\mathrm{s}(\lambda,\alpha,\beta)(2\pi/\kappa-s)}P^\s(\lambda,\alpha,\beta) e_2 h(s) \,\rmd s\
+ \int_0^x \rme^{\eta^\mathrm{s}(\lambda,\alpha,\beta)(x-s)}P^\s(\lambda,\alpha,\beta) e_2 h(s) \, \rmd s,
\end{align*}
and the center and stable spectral projections are of the form
\begin{align*}
P^\mathrm{c}(\lambda,\alpha,\beta) =
\frac{1}{1-\frac{\eta^\mathrm{c}}{\eta^\mathrm{s}}}
\begin{pmatrix} 1 \\ \eta^\mathrm{c} \end{pmatrix}
\begin{pmatrix} 1,-\frac{1}{\eta^\mathrm{s}} \end{pmatrix}, \qquad
P^\mathrm{s}(\lambda,\alpha,\beta) = 
\frac{1}{1-\frac{\eta^\mathrm{c}}{\eta^\mathrm{s}}}
\begin{pmatrix} \frac{1}{\eta^\mathrm{s}} \\ 1 \end{pmatrix}
\begin{pmatrix} -\eta^\mathrm{c},1 \end{pmatrix}.
\end{align*}

To solve the fixed-point problem (\ref{e:fpe1}), it suffices to show that $\mathcal{J}B_3\mathcal{D}$ lies in $L(X^0)$ and has norm strictly less than one. In addition to proving this contraction property, we will derive an expansion of $\mathcal{J}$ as this will help us find expansions of the dispersion curve. We define
\begin{align} \label{def:jc}
[\mathcal{J}^\mathrm{c}_0(\lambda)h](x) := \frac{1}{\omega} \int_0^x \rme^{(\lambda-\bar{g})(x-s)/\omega} h(s) \,\rmd s
\end{align}
and set
\begin{align*}
\mathcal{J}(\lambda,\alpha,\beta) = \mathcal{J}^\mathrm{c}(\lambda,\alpha,\beta) + \mathcal{J}^\mathrm{s}(\lambda,\alpha,\beta), \qquad
\mathcal{J}^\mathrm{c,s}(\lambda,\alpha,\beta) := \frac{1}{\alpha^2\beta^2} P_1 P^\mathrm{c,s} J(\lambda,\alpha,\beta).
\end{align*}
Using the identities
\begin{align}\label{e:p1pce2}
\frac{1}{\alpha^2\beta^2} P_1 P^\mathrm{c}(\lambda,\alpha,\beta) e_2 =
\frac{-1}{\alpha^2\beta^2(\eta^\mathrm{s}+\eta^\mathrm{c})} = 
\frac{1}{\omega} + \rmO(\alpha\beta^2), \qquad
\frac{1}{\alpha^2\beta^2} P_1 P^\mathrm{s}(\lambda,\alpha,\beta) e_2 =
\frac{-1}{\omega} + \rmO(\alpha\beta^2)
\end{align}
and proceeding as in the proof of Lemma~\ref{l:2}, we obtain
\begin{align}
\|\mathcal{J}^\mathrm{c}(\lambda,\alpha,\beta)\|_{L(X^0)} & =
\left\|\frac{1}{\alpha^2\beta^2} P_1 P^\mathrm{c} J\right\|_{L(X^0)}
\leq C \frac{|P_1 P^\mathrm{c} e_2|}{\alpha^2\beta^2}
\leq C
\label{e:jc} \\ \label{e:js}
\|\mathcal{J}^\mathrm{s}(\lambda,\alpha,\beta)\|_{L(X^0)} & =
\left\|\frac{1}{\alpha^2\beta^2} P_1 P^\mathrm{s} J\right\|_{L(X^0)}
\leq \frac{|P_1 P^\mathrm{s} e_2|}{\alpha^2\beta^2} \frac{C}{|\Re\eta^\mathrm{s}|}
\leq \frac{C}{|\Re\eta^\mathrm{s}|}
\leq C\alpha^2\beta^2 \\ \label{e:jcomp}
\|\mathcal{J}^\mathrm{c}(\lambda,\alpha,\beta) - \mathcal{J}^\mathrm{c}_0(\lambda)\|_{L(X^0)} & \leq C \beta^2.
\end{align}
Hence, $\|\mathcal{J}(\lambda,\alpha,\beta)\|_{L(X^0)}\leq C$, and a similar estimate shows that we also have $\|\mathcal{J}_\lambda(\lambda,\alpha,\beta)\|_{L(X^0)}\leq C$, which we will use below. Using the estimate (\ref{e:d}) for $\mathcal{D}$, we conclude that
\begin{align*}
\|\mathcal{J}(\lambda,\alpha,\beta) B_3 \mathcal{D}(\lambda,\alpha)\|_{L(X^0)} \leq C |g_1|_{X^0} C_1 \alpha \leq C\alpha.
\end{align*}
Thus, there is an $\alpha_4>0$ such that (\ref{e:fpe1}) has a unique (up to scalar multiples) nontrivial solution $v\in X^2$ for each $(\lambda,\alpha,\beta)$ with $\lambda\in\Lambda_R$, $0<\alpha\leq\alpha_4$, and $0<\beta\leq1$, and upon setting $a^\mathrm{c}=(1,\eta^\mathrm{c})^t$ this solution is given by
\begin{align*}
v = \left(1+\mathcal{J}(\lambda,\alpha,\beta)B_3\mathcal{D}(\lambda,\alpha)\right)^{-1} \rme^{\eta^\mathrm{c}(\lambda,\alpha,\beta)x}
= \left(1+\mathcal{J}^\mathrm{c}_0(\lambda)B_3\mathcal{D}(\lambda,\alpha)\right)^{-1} \rme^{(\lambda-\bar{g})x/\omega} + \rmO(\beta^2)
= 
\end{align*}
where we used (\ref{e:nu}), (\ref{e:js}), and (\ref{e:jcomp}).

It remains to solve the periodicity condition (\ref{e:1D}) with $a^\mathrm{c}$ as given above. Since $P_1P^\mathrm{c}=1$ on $E^\mathrm{c}$, we can apply $P_1$ to (\ref{e:1D}), which upon using again the estimates (\ref{e:nu}), (\ref{e:js}), and (\ref{e:jcomp}) results in the equivalent equation
\begin{align}
\Delta(\lambda,\alpha,\beta) & :=
\rme^{2\pi\eta^\mathrm{c}(\lambda,\alpha,\beta)/\kappa} - 1
- \left[ \mathcal{J}^\mathrm{c}(\lambda,\alpha,\beta) B_3 \mathcal{D}(\lambda,\alpha) \left(1+\mathcal{J}(\lambda,\alpha,\beta) B_3 \mathcal{D}(\lambda,\alpha)\right)^{-1} \rme^{\eta^\mathrm{c}(\lambda,\alpha,\beta)x} \right](2\pi/\kappa)
\nonumber \\ \label{e:3D} & =
\rme^{2\pi(\lambda+\beta^2-\bar{g})/(\kappa\omega)} - 1
- \left[ \mathcal{J}^\mathrm{c}_0(\lambda) B_3 \mathcal{D}(\lambda,\alpha) \left(1+\mathcal{J}^\mathrm{c}_0(\lambda) B_3 \mathcal{D}(\lambda,\alpha)\right)^{-1} \rme^{(\lambda-\bar{g})x/\omega} \right](2\pi/\kappa) + \rmO(\alpha\beta^2)
\\ \nonumber & = 0
\end{align}
that we need to solve. Comparing (\ref{def:jc}) and (\ref{e:3D}) with the expressions (\ref{def:c}) and (\ref{def:q}) for the solutions of the $\delta=0$ equation, we see that we can write (\ref{e:3D}) as
\begin{align}\label{e:4D}
\Delta(\lambda,\alpha,\beta) = \rme^{2\pi(\lambda+\beta^2-\bar{g})/(\kappa\omega)} - 1
- \left[ \mathcal{Q}(\lambda,\alpha)\mathcal{C}(\lambda,\alpha) \right](2\pi/\kappa) + \rmO(\alpha\beta^2) = 0.
\end{align}
We showed in the proof of Proposition~\ref{p:2} that the equation $\Delta(\lambda,\alpha,0)=0$ with $\beta=0$ has the unique solution $\lambda=\lambda_0^*(\alpha)=\bar{g}+\rmO(\alpha^2)$. Hence, setting $\lambda=\lambda_0^*(\alpha)+\mu$ with $\mu$ near zero and using the estimates (\ref{e:cq}) for $\mathcal{C}$ and $\mathcal{Q}$, we see that (\ref{e:4D}) can be written as
\begin{align*}
\Delta(\lambda_0^*(\alpha)+\mu,\alpha,\beta) & =
\rme^{2\pi(\lambda_0^*(\alpha)+\mu+\beta^2-\bar{g})/(\kappa\omega)} - 1 - \left[ \mathcal{Q}(\lambda_0^*(\alpha)+\mu,\alpha)\mathcal{C}(\lambda_0^*(\alpha)+\mu,\alpha) \right](2\pi/\omega) + \rmO(\alpha\beta^2) \\ & =
\frac{2\pi}{\kappa\omega} (\mu+\beta^2) + \rmO(|\mu+\beta^2|^2) + \rmO(\alpha\mu) + \rmO(\alpha\beta^2) = 0,
\end{align*}
which we can solve uniquely for $\mu=\mu^*(\alpha,\beta)$ with $\mu^*(\alpha,\beta)=-\beta^2(1+\rmO(\alpha))$. In summary, we proved that (\ref{e:3D}) has a solution $(\lambda,\alpha,\beta)\in\Lambda_R\times(0,\alpha_4],\times(0,\beta_4]$ if and only if
\begin{align}\label{e:l*}
\lambda = \lambda^*(\alpha,\beta) \quad\mbox{with}\quad
\lambda^*(\alpha,\beta) = \lambda_0^*(\alpha) - \beta^2 \left(1+\rmO(\alpha)\right).
\end{align}

Finally, for each fixed $\delta \geq 0$, we parameterize the curve $\delta= \alpha^2 \beta^2$ in order to write $\lambda^* = \lambda^*(s,\delta)$. Define
\begin{align} \label{e:params}
\left(\alpha(s,\delta),\beta(s,\delta)\right) = \frac{1}{2}\left( s + \sqrt{s^2 + 4\sqrt{\delta}}, -s + \sqrt{s^2 + 4\sqrt{\delta}} \right)
\end{align}
for $|s|\leq 1$. Then,
\begin{align} \label{e:ab0}
\left( \alpha(s,0), \beta(s,0) \right) = \begin{cases} (s,0) & s \geq 0 \\ (0,|s|) & s \leq 0 \end{cases}, \qquad \text{ and }\qquad  \max_{|s| \leq 1} \left|(\alpha,\beta)(s,\delta) - (\alpha,\beta)(s,0) \right| \leq \delta^{1/4}
\end{align}
Using (\ref{e:params}), write $\lambda^*(s,\delta) := \lambda^*\left(\alpha(s,\delta), \beta(s,\delta) \right)$, and from (\ref{e:l*}) we find
\begin{align} \label{e:lambda_sd}
\lambda^*(s,\delta) &= \lambda^0_*\left(\alpha(s,\delta)\right) - \beta^2(s,\delta)\left(1 + \rmO(\alpha(s,\delta))\right).
\end{align}
For $\delta$ small, we have 
\begin{align*}
\lambda_0^*(\alpha(s,\delta)) - \lambda_0^*(\alpha(s,0)) = \begin{cases}\lambda_0^*(\alpha(s,\delta)) - \lambda_0^*(s) & s \geq 0\\ \lambda_0^*(\alpha(s,\delta)) - \bar{g} & s \leq 0  \end{cases}
\end{align*}
hence
\begin{align} \label{e:d0}
|\lambda_0^*(\alpha(s,\delta)) - \lambda_0^*(\alpha(s,0))| = \rmO(\alpha^2) = \rmO(\delta^{1/2}), \ \ \forall |s| \leq 1.
\end{align}
Thus, using formulation (\ref{e:lambda_sd}) along with estimates in (\ref{e:ab0}) and (\ref{e:d0}) and we find that there exists $s_*,\delta_* >0$ such that
\begin{align*}
\lambda^*(s,\delta) = \begin{cases} \lambda^0_*\left(s\right) + \rmO(\delta^{1/2}) & s \geq 0 \\  \bar{g} - s^2 \left(1 + \rmO(\delta^{1/4})\right) + \rmO(\delta^{1/2})&  s \leq 0  \end{cases}
\end{align*}
uniformly in $|s| \leq s_*$ and $0<\delta<\delta_*$. Hence, we have proved that (\ref{evp:7})-(\ref{evp:8}) has a non-trivial $2\pi/\kappa$-periodic solution if and only if $\lambda = \lambda_*(s,\delta)$ as defined above. 
\end{Proof}



\section{Predictions for spiral waves on bounded disks when $\delta = 0$} \label{sec:abs_spec}

The essential spectrum does not supply stability information for spirals on bounded disks of radius $R$; instead the absolute spectrum is relevant. The absolute spectrum provides limit points to which infinitely many discrete eigenvalues converge to as $R \rightarrow \infty$ \cite{Sandstede:2000ug,Sandstede:2000ut}. For bounded domains, the essential spectrum nevertheless furnishes useful knowledge of the spatial eigenvalue distribution, which is exploited in Theorem~\ref{thm:abs_spec} to predict absolute spectrum locations. As the absolute spectrum and spatial eigenvalues are described in detail in \cite{Fiedler:2000,Sandstede:2000ug,Sandstede:2000ut,Sandstede:2002ht,Kapitula:2013}, we provide only a brief description of these topics. 

For each $\lambda \in \mathbb{C}$, there exist infinitely many spatial eigenvalues $\nu \in \mathbb{C}$ and eigenfunctions that satisfy relation~(\ref{eqn:sp_disp_rel}). Spatial eigenvalues $\nu$ depend continuously on $\lambda$. For $\lambda \gg 1$, we divide $\nu$ into two disjoint groups labeled \enquote{left} for those $\nu$ with $\Re \nu <0$ and \enquote{right} for those $\nu$ with $\Re \nu >0$. Spatial eigenvalues $\nu$ will retain their labels even if their real part changes sign as $\lambda$ varies. The essential spectrum is defined by the set of $\lambda$ for which there exists a purely imaginary spatial eigenvalue. When crossing a positively (negatively) orientated essential spectrum curve from right to left, a right (left) spatial eigenvalue crosses the imaginary axis. The absolute spectrum is defined as values of $\lambda$ for which one right and one left spatial eigenvalue have precisely the same real part: the absolute spectrum must necessarily lie to the left of an essential spectrum branch.

For the scenario $\delta = 0$, there are three essential spectrum curves that influence the possible locations of the absolute spectrum. We denote the essential spectrum curve that is tangent to the imaginary axis by $\Sigma^0_{\text{ess}}$ and the two curves that form the cusp to be $\Sigma^{c,1}_{\text{ess}}$ and $\Sigma^{c,2}_{\text{ess}}$. Next, we denote the region of the $\lambda$-plane enclosed by $\Sigma^{0}_{\text{ess}}$, $\Sigma^{c,1}_{\text{ess}}$, and $\Sigma^{c,2}_{\text{ess}}$ by $\Omega$. Furthermore, we label the spatial eigenvalues that lie on $\rmi \mathbb{R}$ along these essential spectrum curves as $\nu^-_{1}, \ \nu^-_{2},$ and  $\nu^+_{1}$ where the superscripts give the sign of Re $\nu_j$ for $\lambda \gg 1$ and subscripts designate the order the spatial eigenvalues cross the imaginary axis. 

\begin{figure}
    \centering
    \begin{subfigure}{0.45\textwidth}
        \includegraphics[width=\textwidth]{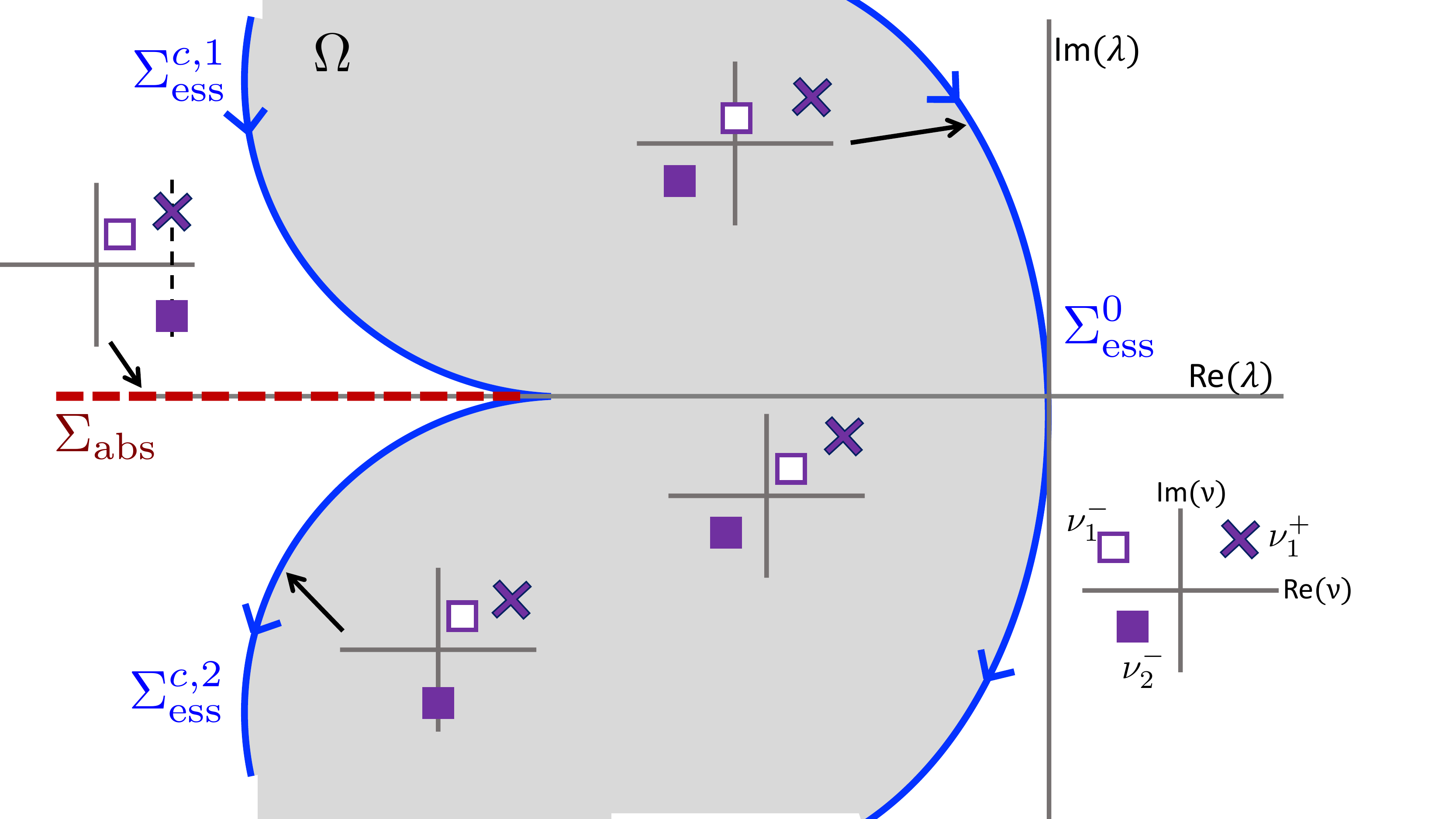}
\caption{Case 1}
        \label{fig:ess_spec_case1}
    \end{subfigure}
    \begin{subfigure}{0.45\textwidth}
        \includegraphics[width=\textwidth]{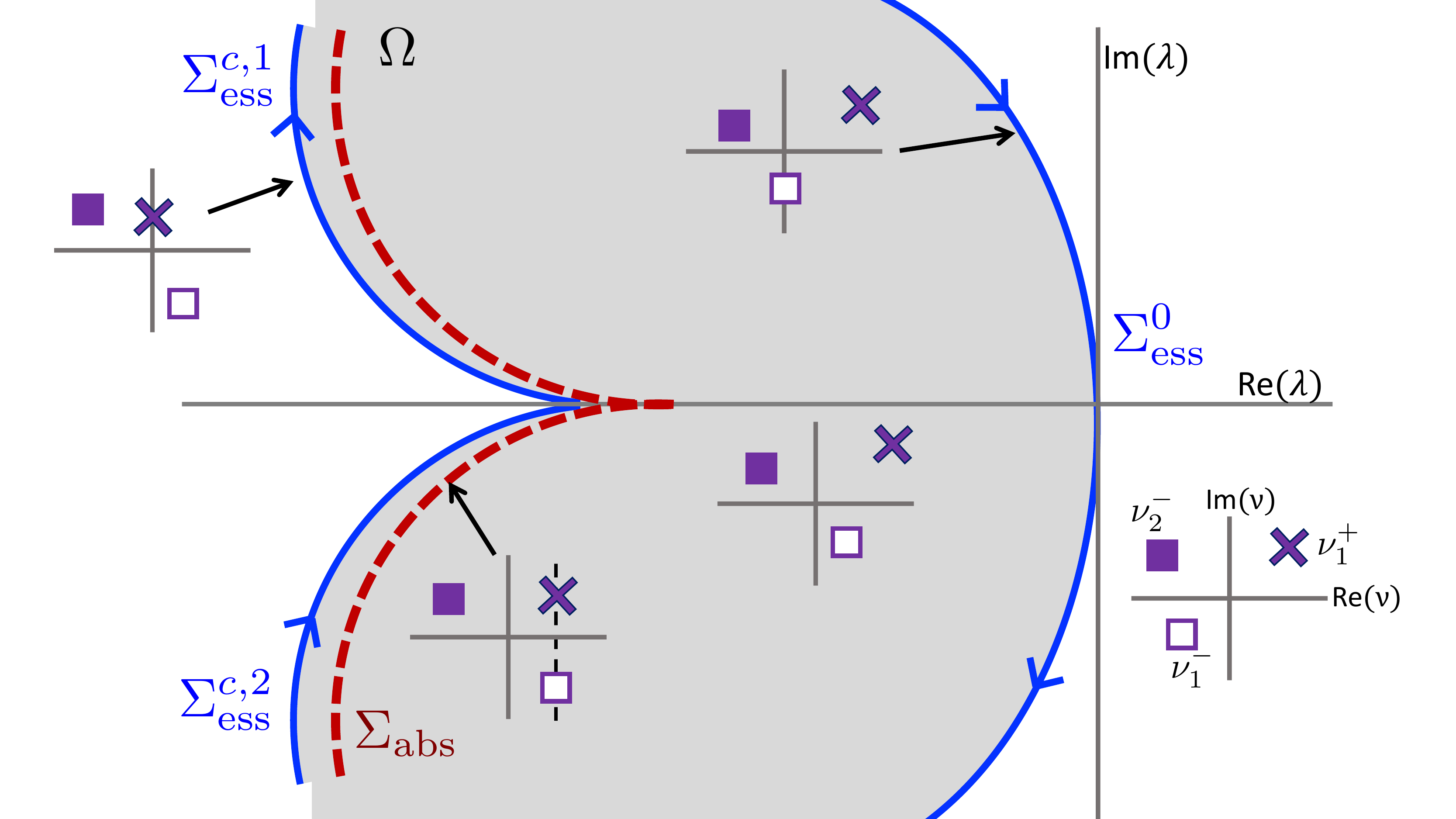}
\caption{Case 2}
        \label{fig:ess_spec_case2}
    \end{subfigure}
    \caption{Schematic illustration for two cases of essential spectrum branches. Arrows on essential spectrum branches indicate the orientation. Insets show the distribution of spatial eigenvalues $\nu(\lambda)$: crosses (squares) indicate spatial eigenvalues $\nu(\lambda)$ that have positive (negative) real part for $\lambda\gg 1$. Shaded region between essential spectrum curves is denoted $\Omega$. (a) Case 1 ($\lambda_3 >0$): Same orientation of essential spectrum branches, leading to the absolute spectrum falling to the left of the cusp point. (b) Case 2 ($\lambda_3 <0$): Opposite orientation of essential spectrum branches, leading to the absolute spectrum to the right of the cusp point.}\label{fig:figs34}
\end{figure}

Two separate cases emerge for the orientation of essential spectrum branches that compose the cusp structure at $\lambda_0$. The curve orientation is dictated by the sign of the imaginary $\alpha^3$ term in $\lambda_*^0(\alpha)$. A positive coefficient (Case 1) leads to $\Sigma^{c,1,2}_{\text{ess}}$ having the same orientation as $\Sigma^0_{\text{ess}}$, and a negative coefficient (Case 2) produces curves of opposite orientation (see Figure~\ref{fig:figs34}). In both cases, $\Sigma^0_{\text{ess}}$ is negatively oriented and crossing this curve from right to left results in the  spatial eigenvalue $\nu^-_1$ crossing the imaginary axis (Figure~\ref{fig:figs34}). Theorem~\ref{thm:abs_spec} predicts the locations of the absolute spectrum for each case. 

\begin{Theorem}\label{thm:abs_spec} 
Consider the essential spectrum $\lambda_*^0(\alpha)$ for $\delta = 0$. For $\lambda_2 <0$, predictions of absolute spectrum locations are as follows.
	\begin{enumerate}
		\item Case 1 ($\lambda_3 > 0$): If Re $\nu^-_1(\lambda) < $ Re $\nu^+_1(\lambda)$ for all $\lambda \in \Omega$, then there are no absolute spectra curves within $\Omega$. Moreover, locally the absolute spectrum is given by the set
\begin{align*}
 \Sigma_{\text{abs}} = \left\{ \lambda \in \mathbb{C} : \lambda = \lambda_0 - r^2, \ \ 0< r< r_0\right\}
\end{align*}
for some small real $r_0>0$.

		\item Case 2 ( $\lambda_3 < 0$): The set $\Omega$ contains a curve of absolute spectrum
	\end{enumerate}
\end{Theorem}

The proof of Theorem~\ref{thm:abs_spec} follows from arguments of the orientation of essential spectrum curves and crossings of spatial eigenvalues. If $\lambda_3 = 0$, the result holds by replacing $\lambda_3$  with the leading order non-zero imaginary term in $\lambda_*^0(\alpha)$.

\begin{Proof}
 \textbf{Case 1:} First, the orientation of essential spectrum curves and assumption on the spatial eigenvalues ensures there is no absolute spectrum curve within $\Omega$ (see schematic in Figure~\ref{fig:figs34}a). An absolute spectrum curve will be given by the set of $\lambda$ for which left and right spatial eigenvalues have the same positive real part. 

We use $\lambda_*^0(\alpha)$ to solve locally for $\alpha = \alpha(\lambda)$ to find the spatial eigenvalues given by $\nu = \rmi \gamma = \rmi/\alpha$. We have $\lambda_*^0(\alpha) = \lambda_0 -  \lambda_2 \alpha^2 + \rmi  \lambda_3 \alpha^3$ which is written such that $\lambda_2, \lambda_3 >0$. After shifting the curve by $\lambda_0$, we want to find solutions $\alpha(\lambda)$ of 
\begin{align*}
\mathcal{F}(\alpha,\lambda) = -\lambda_2 \alpha^2 + \rmi \lambda_3 \alpha^3 - \lambda = 0
\end{align*}
and determine for which values of $\lambda$ two spatial eigenvalues have the same real part. Using a Newton's polygon, set $\alpha = \lambda^{1/2}w$ which transforms the equation to
\begin{align*}
-\lambda_2 w^2+ \rmi \lambda_3 \lambda^{1/2}w^3 - 1 = 0.
\end{align*}
Near $\lambda = 0$, $w^2 = -1/\sqrt{\lambda_2}$ and there are two unique roots $\alpha_{1,2} = \lambda^{1/2}w_{1,2}(\lambda^{1/2})$ where $w_{1,2}$ are analytic in $\lambda^{1/2}$ and $w_{1,2}(0) = \pm \rmi/\sqrt{\lambda_2}$. The final root is found by setting $\alpha = \lambda^0 w(\lambda^0)$, resulting in $\alpha_3(0) = -\rmi \lambda_2/\lambda_3$. Mapping these $\alpha$-roots to $\nu = \rmi /\alpha$ gives
\begin{align*}
\nu_1(\lambda) &= \frac{\rmi}{\lambda^{1/2} v_1(\lambda^{1/2})}, \ \ \nu_1(0) = \frac{\sqrt{\lambda_2}}{\lambda}\\
\nu_2(\lambda) &= \frac{\rmi}{\lambda^{1/2} v_2(\lambda^{1/2})}, \ \ \nu_2(0) = -\frac{\sqrt{\lambda_2}}{\lambda}\\
\nu_3 &= -\frac{\lambda_3}{\lambda_2}.
\end{align*}
The two relevant spatial eigenvalues are $\nu_1$ and $\nu_2$. Using the leading order in $\lambda$ terms, set $\Re \nu_1 =  \Re \nu_2$  and let $\lambda = r^2 \rme^{2\rmi \theta}$. Then, solving
\begin{align*}
\frac{\sqrt{\lambda_2}}{r \cos(\theta)} = -\frac{\sqrt{\lambda_2}}{r \cos(\theta)}
\end{align*}
gives $\theta = \pi/2$ or $\lambda = r^2 \rme^{\rmi \pi} = -r^2$. Locally, near $\lambda_0$, there exists a small $r_0 >0$ such that an absolute spectrum curve is given by 
\begin{align*}
 \Sigma_{\text{abs}} = \left\{ \lambda \in \mathbb{C} : \lambda = \lambda_0 - r^2, \ \ 0 < r< r_0 \right\}.
\end{align*}

\textbf{Case 2:} The results for Case 2 follow from the orientations of the essential spectrum curves. Upon crossing $\Sigma_{\text{ess}}^0$ from right to left a left spatial eigenvalue crosses the imaginary axis. Crossing $\Sigma^{c,1,2}_{\text{ess}}$ from right to left results in a right spatial eigenvalue crossing the imaginary axis. Since spatial eigenvalues are analytic in $\lambda$, there must be values of $\lambda \in \Omega$ for which the real parts of the right and left spatial eigenvalues coincide to give an absolute spectrum curve.
\end{Proof}

In Case 1, we prove the existence of an absolute spectrum branch locally near $\lambda_0$ which extends leftward in the complex plane. Theorem~\ref{thm:abs_spec} guarantees the existence of a spectral curve within $\Omega$ for Case 2, signifying that these models can undergo an absolute instability if $\Omega$ reaches into the positive half plane.

\section{Application to Barkley and Karma models} \label{sec:barkley_karma}

In this section, the results of Theorems~\ref{t:1}-\ref{thm:abs_spec} are applied to planar spiral waves formed in the Barkley and Karma models, which are well-studied nonlinear reaction-diffusion systems \cite{Barkley:1994vk,Barkley:1995fv,Wheeler:2006eu,Sandstede:2006ei,Marcotte:2015ex,Marcotte:2016hn,Dodson:2019ku,Karma:1993tq,Karma:1994gb}. As we will see, these two models also furnish examples of the two absolute spectrum cases.

\subsection{Barkley Model}
The Barkley model is given by
\begin{align*}
u_t &= \ \Delta u + \frac{1}{\epsilon} u (1-u) \left( u - \frac{v + b}{a} \right)\\
v_t &= \ \delta \Delta v + u - v, 
\end{align*}
where we fix parameters $a = 0.7$, $b = 0.001$, and $\epsilon = 0.02$. For these parameters, there is numerical evidence that a stable spiral wave exists \cite{Barkley:1992ia}. The effects of removing diffusion are studied by modifying the diffusion coefficient $\delta$ for the slowly diffusing species $v$ in the interval $[0, 0.2]$. The angular frequency changes with $\delta$ and falls between $\omega_0 = 1.87 \ (\delta = 0.2)$ and $\omega_0 = 2.09 \ (\delta = 0)$. Note that $\bar{g} = -1$.

\begin{figure}
\centering
 \includegraphics[width=0.5\linewidth]{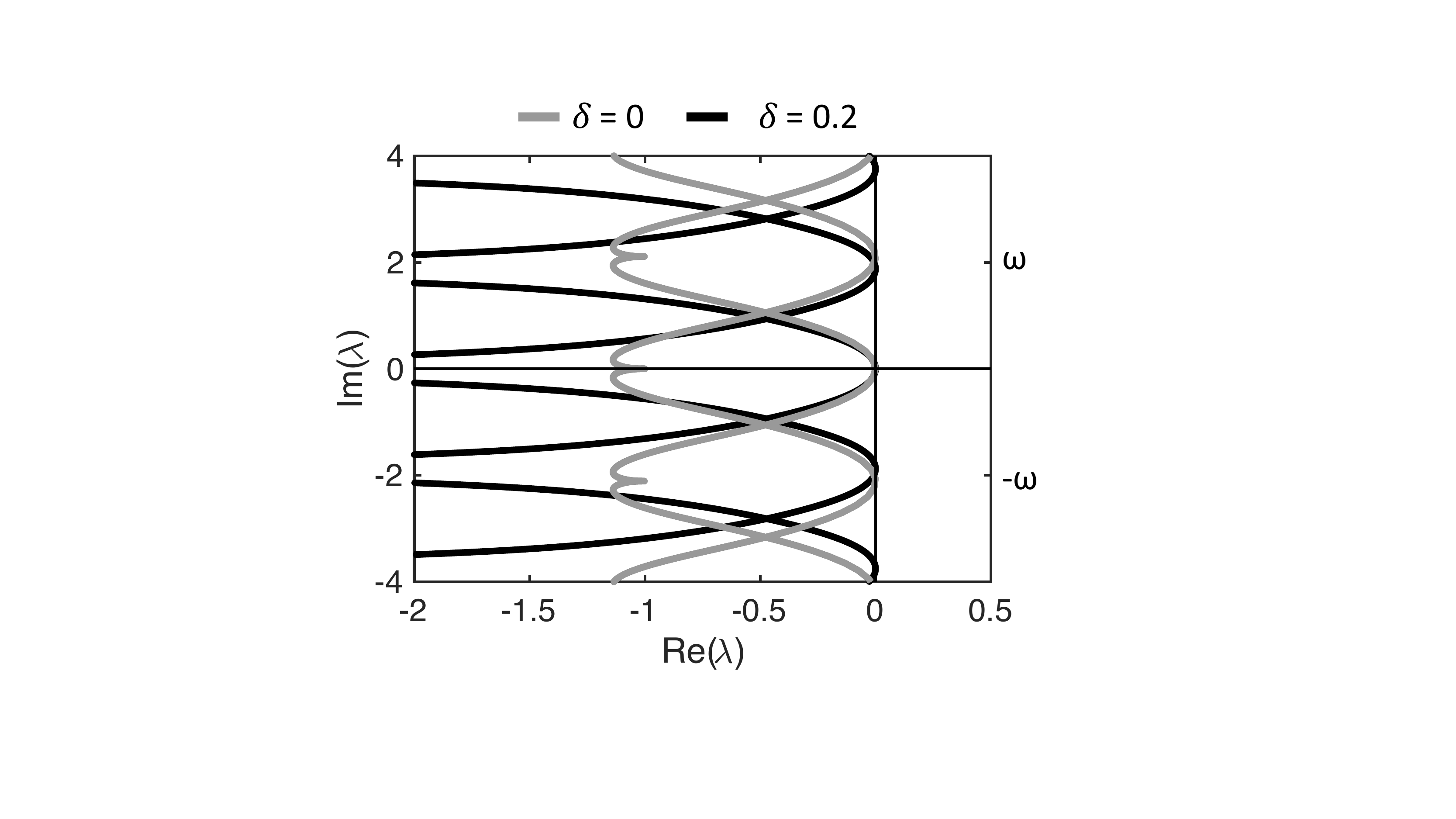}
\footnotesize \caption{Essential spectra for the Barkley model for $\delta = 0$ and $\delta >0$ } \label{fig:bd}
\end{figure}

Essential spectrum curves are calculated using numerical continuation methods on system~(\ref{eqn:2comp_rxn_diff_ess_spec}) as described in \cite{Rademacher:2007uh}. Computations are performed with Matlab on a 256 grid point periodic domain with Fourier spectral methods used for spatial derivatives. Below, these numerically computed spectral curves will be compared to the analytical predictions.

Figure~\ref{fig:bd} shows essential spectrum curves of the spiral wave calculated for $\delta = 0$ and $\delta =0.2$. Vertically periodicity of the branches arises due to the Floquet symmetry of the spiral eigenfunctions. As noted above and as indicated in the theorems, the curves look very different for these two cases. For $\delta=0.2$, the curves are unbounded. In the case $\delta = 0$, adjacent branches meet at the limit points $\lambda_0$ and form cusps. The limit points are predicted to be $\lambda_0 = -1 + \rmi \omega_0 n$.

\begin{figure}
\centering
 \includegraphics[width=0.8\linewidth]{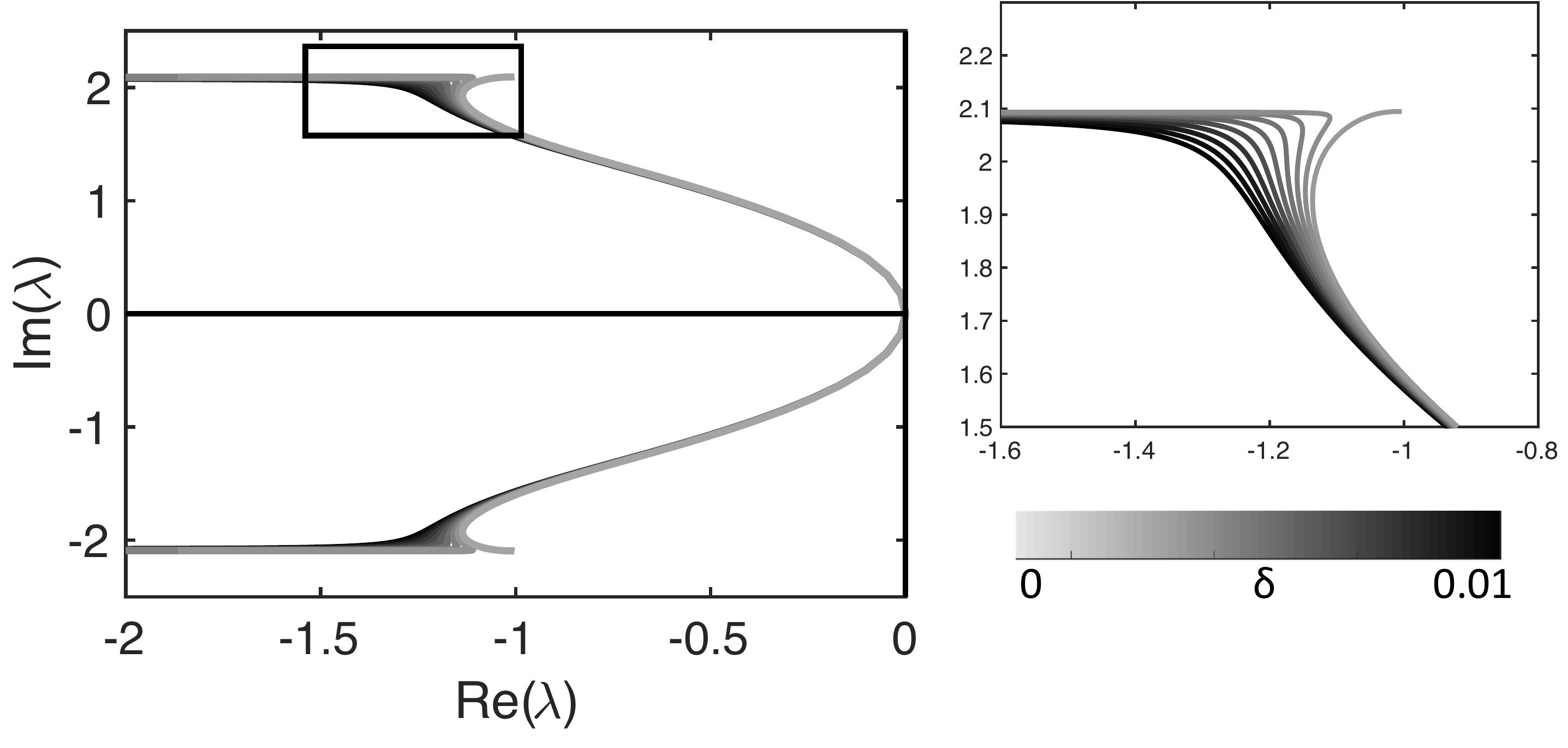}
\footnotesize \caption{Essential spectra for the Barkley model under changing diffusion. Inset at top right focuses on area near $\lambda_0$. } \label{fig:fig5}
\end{figure}

Numerically computed curves for decreasing $\delta$ are shown in Figure~\ref{fig:fig5} with the inset focusing on the region near $\lambda_0$. For small, non-zero values of $\delta$, curves begin to turn toward the $\lambda_0$ limit point before diverging.

We verify the form of $\lambda_*^0(\alpha)$ computed in Theorem~\ref{t:1} by analyzing the convergence of numerically computed spectral curves; denote these curves by $\lambda_c$.  From $\lambda_*^0(\alpha)$, we expect $\Re\left(|\lambda_c - \lambda_0|\right)=\rmO(\alpha^2)$ and $\Im\left( |\lambda_c - \lambda_0| \right)=\rmO(\alpha^3)$. Log-log plots of the convergence are shown in Figure~\ref{fig:fig6} and show the $\rmO(\alpha^2)$ dependence of the real terms. Note that $\Im\left( |\lambda_c - \lambda_0| \right)=\rmO(\alpha^5)$. This leads to the following observation.

\begin{figure}
\centering
 \includegraphics[width=0.4\linewidth]{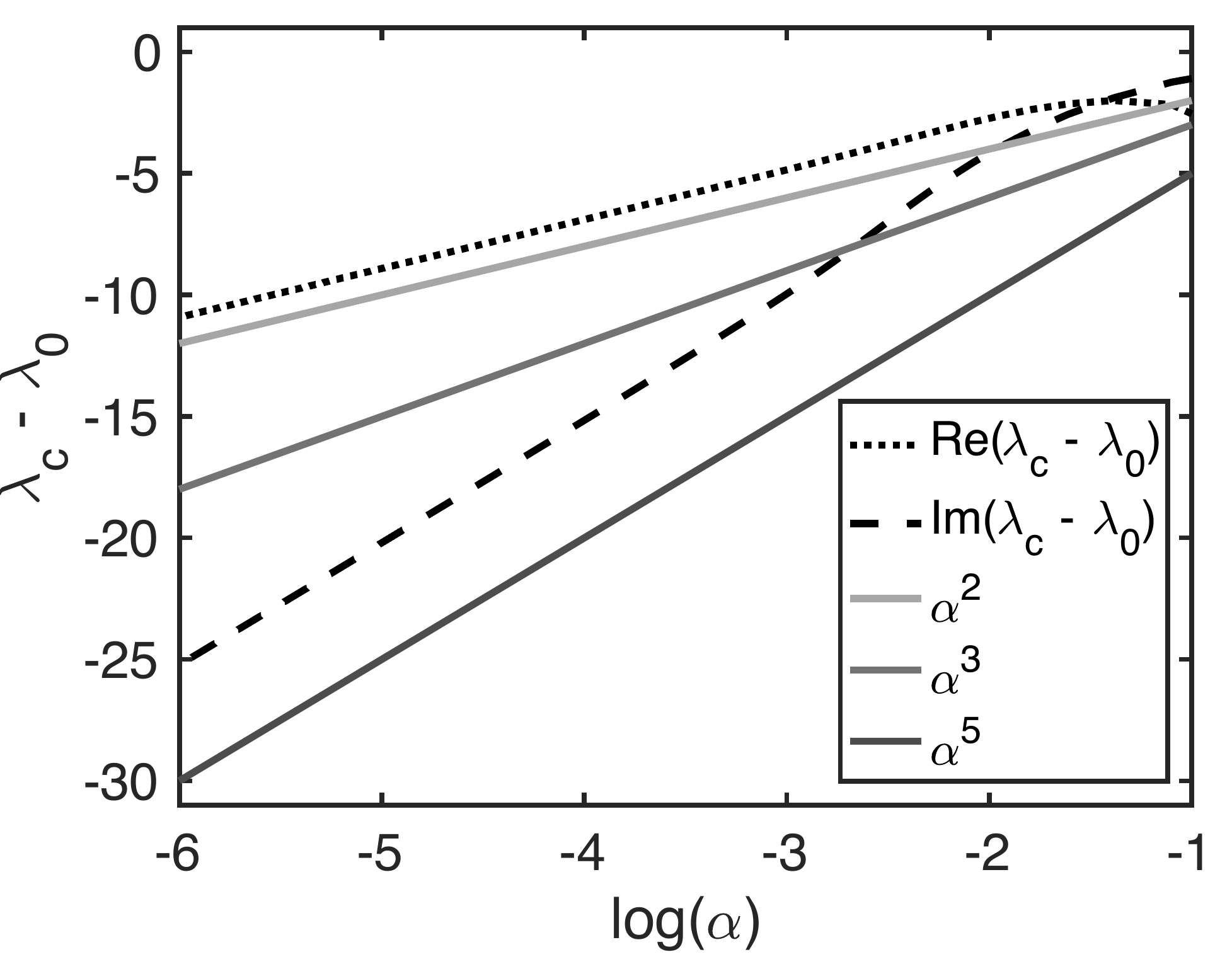}
\footnotesize \caption{ Convergence of real and imaginary parts of the essential spectrum computed with continuation in the Barkley model for $\delta = 0$.} \label{fig:fig8}
\end{figure}

\begin{Corollary} \label{cor:gu_constant}
If the linearized term $g_1(x)$ is a constant, then $\lambda_3 = 0$ in $\lambda_*^0(\alpha)$.
\end{Corollary}
This corollary is immediately proved by integrating $\lambda_3$ and using the periodicity of $f_2(x)$. The linearization $g_1(x) = 1$ in the Barkley model, and Corollary \ref{cor:gu_constant} therefore applies.

Numerical computation of absolute spectra following methods in \cite{Rademacher:2007uh} for $\delta = 0$ and $\delta =0.2$ are shown in Figure~\ref{fig:fig1}. There is an exceptional difference between the two: the Y-shaped branches present for $\delta = 0.2$ are replaced by short segments that lie to the left of the essential spectrum branches for $\delta =0$.

Under the nonlinearities and parameters used here for the Barkley model, the orientation of the essential spectrum curves means the system falls into absolute spectrum case 1 described in Theorem~\ref{thm:abs_spec}. The numerical computations confirm the existence of an absolute spectrum branch that emerges leftward from $\lambda_0$ (Figure~\ref{fig:fig1}a). Moreover, no absolute spectrum is numerically found within $\Omega$.

\subsection{Karma Model}

\begin{figure}
\centering
 \includegraphics[width=0.9\linewidth]{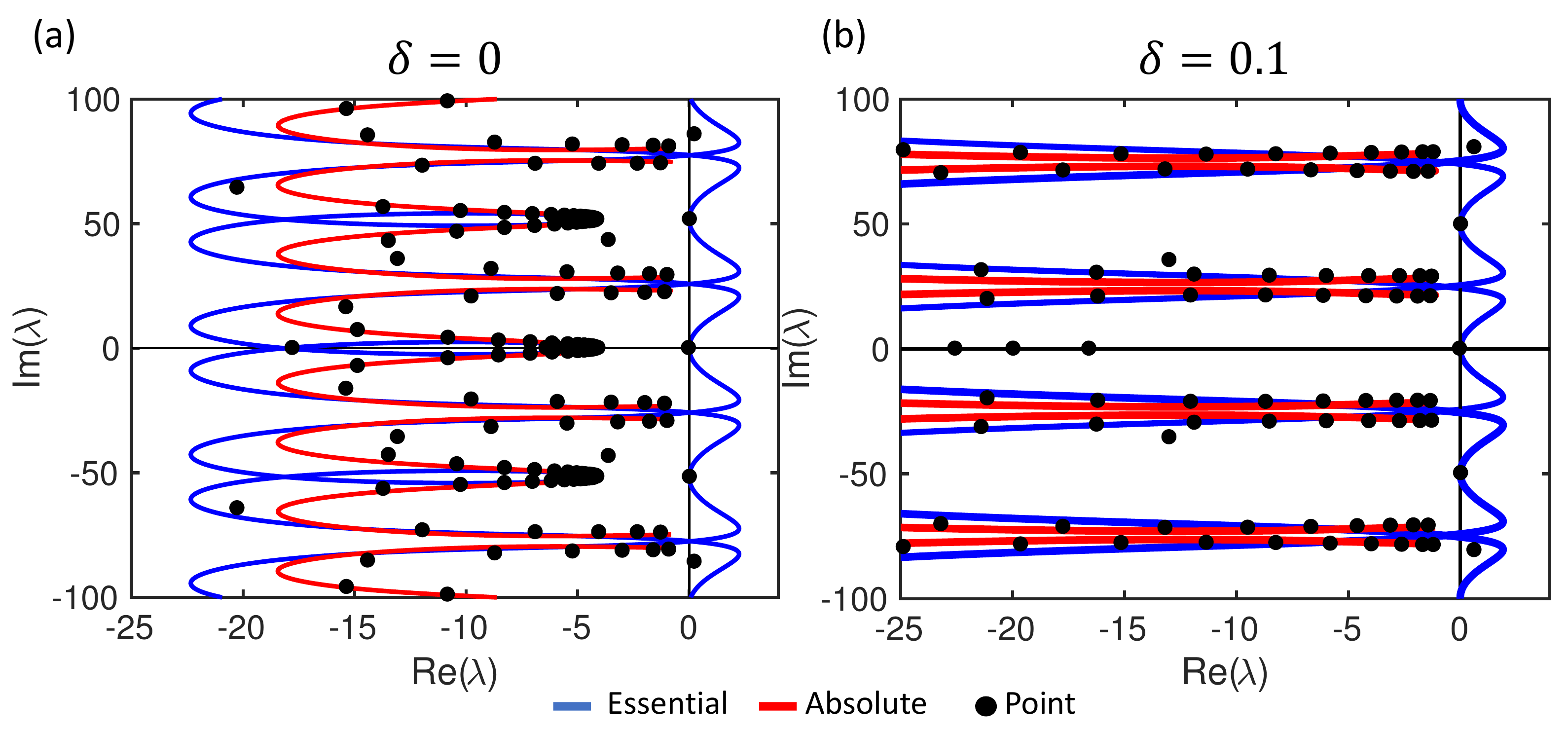}
\footnotesize \caption{Shown are the essential spectra in the Karma model for $\delta=0$ and $\delta=0.1$. When $\delta = 0$, cusps form in the essential spectrum at $\lambda_0 = -4 + \rmi \omega_0 n$. } \label{fig:fig6}
\end{figure}

The Karma model is given by
\begin{align*}
u_t = & \ 1.1 \Delta u + 400 \left( -u + \left( 1.5414 - v^4\right) \left(1 - \tanh(u-3) \right) \frac{u^2}{2} \right)\\
v_t = & \ \delta \Delta v + 4 \left( \frac{1}{1 - \rme^{-\mu_K}} \theta_s(u-1) - v \right)
\end{align*}
where $\mu_K = 1.2$. As in \cite{Allexandre:2004cd,Marcotte:2015ex,Marcotte:2016hn,Dodson:2019ku}, we consider a smoothed version of the Heaviside function given by $\theta_s(u)=(1+\tanh(su))/2$ for $s = 4$. Alternans are numerically observed for this set of parameter values \cite{Marcotte:2015ex,Marcotte:2016hn,Dodson:2019ku} and spirals have angular frequencies of $\omega_0 = 49.81$ ($\delta = 0.1$) and $\omega_0 = 51.66$ ($\delta = 0$). Again, there are distinct differences between the essential and absolute spectra for $\delta >0$ and $\delta = 0$ in Figure~\ref{fig:fig6}. Here, $\bar{g} = -4$.

\begin{figure}
\centering
 \includegraphics[width=0.75\linewidth]{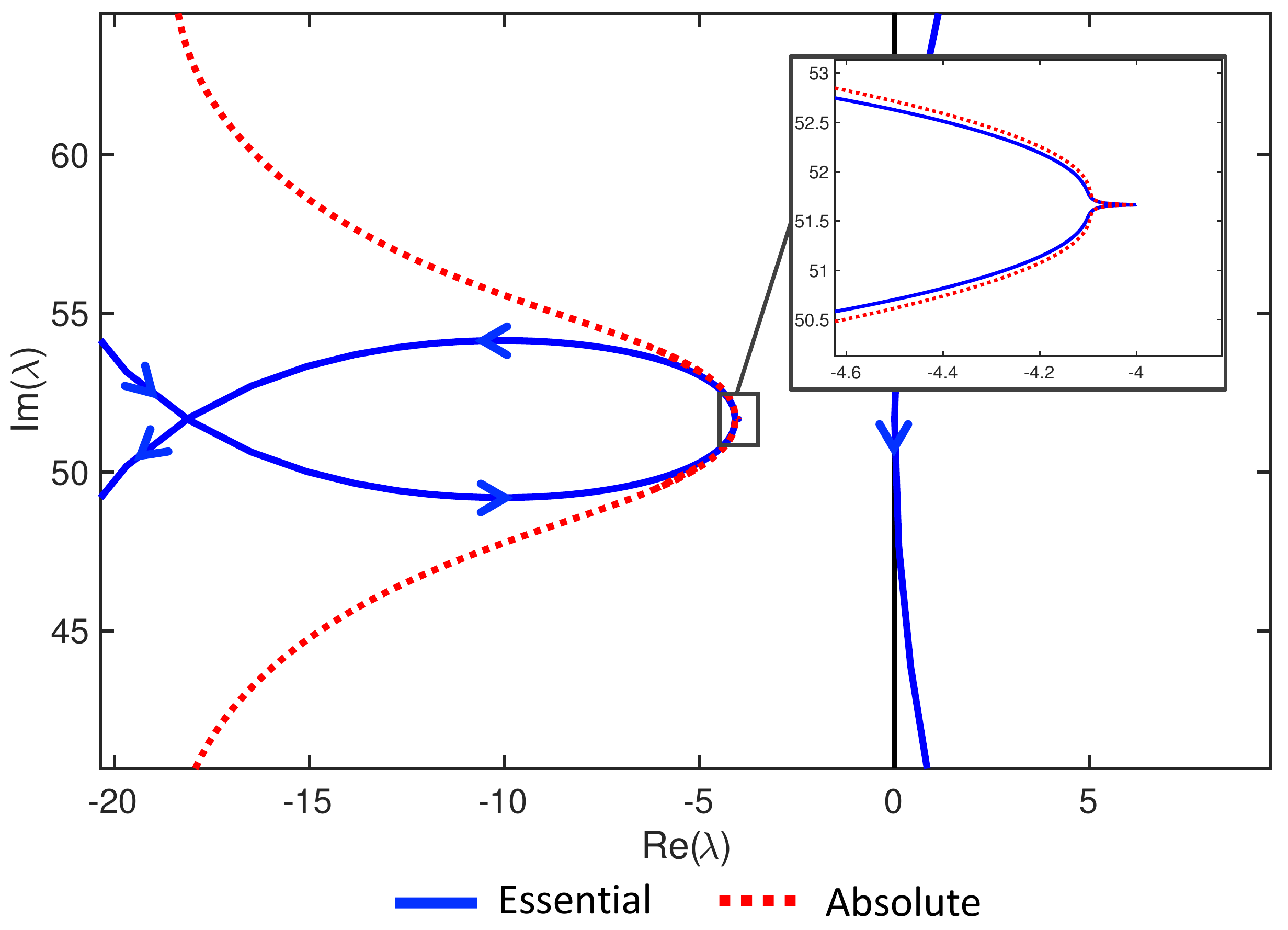}
\footnotesize \caption{Formation of cusps in Karma model for $\delta = 0$. Loops in the essential spectrum change the local orientation (indicated by arrows) at the cusp points. Spectra computed by numerical continuation. } \label{fig:fig7}
\end{figure}

From Theorem~\ref{t:1} the cusps are predicted to emerge from $\lambda_0 = -4 + \rmi\omega_0 n$, which is confirmed numerically (Figure~\ref{fig:fig7}). For $\delta = 0$, essential spectrum curves in the Karma model form loops that reverse the branch orientations before forming the cusp point, and a computation shows that $\lambda_3 <0$. Thus, the Karma model falls into absolute spectrum case 2, and Theorem~\ref{thm:abs_spec} predicts absolute spectrum to the right of $\lambda_0$ in region $\Omega$. The numerical computation of the absolute spectrum shown in Figure~\ref{fig:fig7} confirms this prediction.

Furthermore, the existence of leading absolute spectrum branches in $\Omega$ grants the possibility of the system undergoing an absolute instability. In fact, regardless of $\delta$ these branches do destabilize when the system parameter $\mu_k$ is increased above 1.4. 

\section{Discussion}

\paragraph{\textbf{Summary}}
The central question we pursued is how the continuous spectra of planar spiral waves behave in the limit of vanishing diffusion in one of the components of a reaction-diffusion system. We found that these spectra behave discontinuously near cusp points of the essential spectrum of the zero-diffusion limit. We also priovide explicit expansions of these spectra in the spatial Floquet exponent $1/\alpha$ and the small diffusion coefficient $\delta$ that hold for $0\leq\alpha\ll1$ and $0\leq\delta\ll1$ and whose coefficients are determined by the diffusionless variable. We verified these predictions numerically in the Barkley and Karma models.

We also discussed the absolute spectra of spiral waves on bounded disks for $\delta=0$. The sign of $\lambda_3$ from $\lambda_*^0(\alpha)$ divides models into the two distinct cases described in Theorem~\ref{thm:abs_spec}. Absolute spectra lie strictly outside of the region $\Omega$ if $\lambda_3 >0$. If instead $\lambda_3<0$, then there is absolute spectrum inside $\Omega$. These two cases occur in the Barkley and Karma models as demonstrated in \S\ref{sec:barkley_karma}.

\begin{figure}
\centering
 \includegraphics[width=0.4\linewidth]{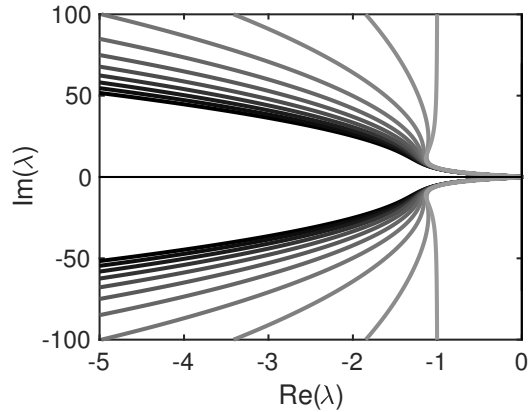}
\footnotesize \caption{Essential spectrum of wave trains in the Barkley model limits to a vertical branch under $\delta \rightarrow 0$.} \label{fig:fig10}
\end{figure}

\paragraph{\textbf{Wave-train spectra}}
As mentioned in the introduction, the essential spectra of the asymptotic wave trains depend smoothly on $\delta$ in the zero-diffusion limit \cite{Rademacher:2004ue}. Interestingly, cusp points are not observed for the essential spectra of these wave trains when $\delta=0$. Figure~\ref{fig:fig10} shows essential spectra for the wave train in the Barkley model for decreasing values of $\delta$. The limiting behavior of the wave train essential spectrum can be obtained directly from the spectral relationship (\ref{eqn:sp_disp_rel}) and is given by
\begin{align*}
\lambda^0_{\infty}(\alpha) = \lambda_0 + \lambda_2\alpha^2 +\rmi \lambda_3\alpha^3 + \omega \left( \frac{\rmi}{\alpha}\right) + \rmo(\alpha^4).
\end{align*} 
Only the imaginary parts of the spectra differ between $\lambda^0_{\infty}(\alpha)$ and $\lambda^0_*(\alpha)$. In the limit $\alpha \rightarrow 0$, we have $\Re \lambda_*^0(\alpha) \rightarrow \bar{g}$, implying $\Re \lambda^0_{\infty}(\alpha) \rightarrow \bar{g}$, but $\Im \lambda^0_{\infty}(\alpha)$ becomes infinite due to the addition of the imaginary $1/\alpha$ term.

For the spiral wave at $\delta = 0$, the limit points $\lambda_0$ arise because the $v$-equation decouples at small values of $\alpha$. To further understand the origins of the differences between wave-train and spiral spectra, consider the linear dispersion relation for the wave train in the co-moving frame 
\begin{align*}
 D \left( \partial_{x} + \nu\right)^2 \bar{V} + \omega\left(\partial_{x} + \nu\right)^2 \bar{V} + F_U(U_{\infty}(x))\bar{V} - \lambda_{\infty} \bar{V} = 0.
\end{align*}
For a two-variable system as studied above with $\delta = 0$ and $\nu = \rmi \gamma$, the $v$-equation takes the form
\begin{align*}
\lambda_{\infty} v = \omega \left( \partial_{x} + \rmi \gamma \right)v + g_1(x) u + \bar{g} v.
\end{align*}
Here, the $v$-equation still retains direct input from $\gamma$, unlike the corresponding equation for the spiral (\ref{e:5}). The $\gamma$ provides an unbounded term for which $\lambda_{\infty}$ must compensate in order to retain bounded eigenfunctions.

\paragraph{\textbf{Outlook}}

There are a number of questions that we did not address. One such question is whether persistence of planar spiral waves that are assumed to exist for $\delta=0$ into the region $0<\delta\ll1$ could be established using singular perturbation theory. Similarly, we discussed only the essential spectra of planar spiral waves in the vanishing-diffusion limit, and the behavior of point spectrum was not addressed here. Moreover, we focused on two-component system, and it is natural to ask whether our results hold for general $n$-component systems, potentially with several non-diffusing variables: we expect that each non-diffusing variable creates its own cusps but have not attempted to prove this. Finally, while Theorem~\ref{t:1} remains true if we allow $\bar{g}$ to be a $2\pi/\kappa$-periodic function, we did not prove the same for Theorem~\ref{t:2}, though we expect this to hold as well.

On bounded disks, the spectrum consists entirely of discrete eigenvalues, which fall into disjoint sets that either align along curves of the absolute spectrum, are members of the extended point spectrum, or arise from the imposed boundary conditions \cite{Sandstede:2000ug,Sandstede:2000ut}. Theorem~\ref{thm:abs_spec} predicts the locations of absolute spectrum depending on the sign of $\lambda_3$, but the effect of vanishing diffusion on boundary and extended point spectrum remains an open question.

We leave the reader with a word of caution: the essential spectrum may be significantly different in the two cases of small positive and zero diffusion $\delta$. This fact has potentially significant consequences for the stability of patterns formed in ion channel models, and care should be taken during computations.

\begin{Acknowledgment}
Dodson was partially supported by the NSF through DMS-1644760.
Sandstede was partially supported by the NSF through grants DMS-1714429 and DMS-2106566. 
\end{Acknowledgment}

\bibliographystyle{abbrvnat}
\bibliography{thesis_ref.bib}

\end{document}